\documentclass[%
 reprint,
superscriptaddress,
preprintnumbers,
 amsmath,amssymb,
 aps,
pre,
]{revtex4-2}

\usepackage{graphicx}
\usepackage{dcolumn}
\usepackage{bm}
\usepackage{hyperref}
\usepackage[mathlines]{lineno}

\usepackage{import}
\usepackage{enumitem} 
\usepackage{cases}
\usepackage{empheq}
\usepackage{siunitx}

\NewDocumentCommand{\HOI}{>{\SplitArgument{2}{}}m}{\HOIaux#1}
\NewDocumentCommand{\HOIaux}{mmm}{$\overleftarrow{#1#2}_{#3}$}

\NewDocumentCommand{\HOIdouble}{>{\SplitArgument{2}{}}m}{\HOIdoubleaux#1}
\NewDocumentCommand{\HOIdoubleaux}{mmm}{$\overleftrightarrow{#1#2}_{#3}$}

\begin{document}

\preprint{APS/123-QED}

\title{Modification speed alters stability of ecological higher-order interaction networks}

\author{Thomas \surname{Van Giel}}
\email{thomas.vangiel@ugent.be}
\author{Aisling J. Daly}%

\affiliation{%
BionamiX, Department of Data Analysis and Mathematical Modelling, Ghent University, 9000 Ghent, Belgium
}%
\affiliation{%
KERMIT, Department of Data Analysis and Mathematical Modelling, Ghent University, 9000 Ghent, Belgium
}%
\author{Jan M. Baetens}
\affiliation{%
BionamiX, Department of Data Analysis and Mathematical Modelling, Ghent University, 9000 Ghent, Belgium
}%

\author{Bernard \surname{De Baets}}

\affiliation{%
KERMIT, Department of Data Analysis and Mathematical Modelling, Ghent University, 9000 Ghent, Belgium
}%

\date{\today}

\begin{abstract}
  Higher-order interactions (HOIs) have the potential to greatly increase our understanding of ecological interaction networks beyond what is possible with established models that usually consider only pairwise interactions between organisms. While equilibrium values of such HOI-based models have been studied, the dynamics of these models and the stability of their equilibria remain underexplored. Here we present a novel investigation on the effect of the onset speed of a higher-order interaction. In particular, we study the stability of the equilibrium of all configurations of a three-species interaction network, including transitive as well as intransitive ones. We show that the HOI onset speed has a dramatic effect on the evolution and stability of the ecological network, with significant structural changes compared to commonly used HOI extensions or pairwise networks. Changes in the HOI onset speed from fast to slow can reverse the stability of the interaction network. The evolution of the system also affects the equilibrium that will be reached, influenced by the HOI onset speed. This implies that the HOI onset speed is an important determinant in the dynamics of ecological systems, and including it in models of ecological networks can improve our understanding thereof.
\end{abstract}

\maketitle

Describing species interactions is at the core of the fields of ecology and ecological modelling. They are crucial in understanding ecology due to their role in shaping community dynamics and responses to environmental changes \cite{akessonImportanceSpeciesInteractions2021}. Historically, species interaction models have typically considered only pairwise interactions. 
Various aspects of these pairwise models have been 
investigated, including transient dynamics~\cite{chenTransientDynamicsFood2001, colonBifurcationAnalysisAgentbased2015}, equilibrium stability~\cite{buninEcologicalCommunitiesLotkaVolterra2017, lotkaAnalyticalNoteCertain1920, reichenbachMobilityPromotesJeopardizes2007, suweisEmergenceStructuralDynamical2013, kesslerGeneralizedModelIsland2015} and ergodicity~\cite{arneodoStrangeAttractorsVolterra1982,vanoChaosLowdimensionalLotka2006}. However, this pairwise framing is clearly a significant simplification of the complexity of a real ecosystem. Therefore, there has been a recent trend in ecological modelling to go beyond pairwise interactions and consider higher-order interactions (HOIs)~\cite{levinePairwiseMechanismsSpecies2017,sanchezDefiningHigherorderInteractions2019, gibbsCoexistenceDiverseCommunities2022} and interaction modification~\cite{terryInteractionModificationsLead2019}. 

Higher-order interactions involve more than two species, and thus cannot be measured using pairwise experiments~\cite{kleinhesselinkDetectingInterpretingHigherorder2022}. These interactions can occur, for example, when a pairwise interaction between two species $A$ and $B$ is modified by the presence of a third species $C$. This modification can be either positive (strengthening the original interaction) or negative (weakening or reversing the original interaction). It can also be asymmetric or symmetric, which means that the third species can either modify only the effect $A$ has on $B$ (asymmetric), or only the effect $B$ has on $A$ (asymmetric), or both (symmetric). Theoretical research has shown that HOIs can have stabilising effects on ecosystems~\cite{grilliHigherorderInteractionsStabilize2017}, and either positive or negative effects on coexistence~\cite{baireyHighorderSpeciesInteractions2016a, terryInteractionModificationsLead2019,xiaoHigherorderInteractionsMitigate2020a}. Thus, theoretical research suggests that HOIs may play a significant role in the evolution and stabilization of ecological systems.

Recent work has focused on detecting HOIs from ecological data~\cite{kleinhesselinkDetectingInterpretingHigherorder2022,rajSimplicialStructuresEcological2022,terryIdentifyingImportantInteraction2020} and other domains~\cite{musciottoDetectingInformativeHigherorder2021}, as well as systems with HOIs in equilibrium~\cite{singhHigherOrderInteractions2021,mayfieldHigherorderInteractionsCapture2017}. Analytical methods have been used for such complex systems like the cavity method~\cite{gibbsCoexistenceDiverseCommunities2022}. Although the dynamics of HOI systems have been studied~\cite{chatterjeeControllingSpeciesDensities2022}, the literature on this topic is relatively scarce, with no studies to our knowledge on stability changes of the equilibrium.

In the real world, HOIs are being discovered more and more frequently in ecological systems~\cite{barbosaExperimentalEvidenceHidden2023,liDirectNeighbourhoodEffects2021,mayfieldHigherorderInteractionsCapture2017,vanveenPlantmodifiedTrophicInteractions2015,morinCompetitionAquaticInsects1988, golubskiEcologicalNetworksEdge2016}. HOIs can manifest in multiple forms: species excreting chemicals that change other species' behaviour~\cite{kumarNaturalHistorydrivenPlantmediated2014}, species changing their foraging behaviour in the presence of a third species~\cite{wissingerIntraguildPredationCompetition1993} or changes in species' competitiveness due to parasitism~\cite{clayEffectsInsectHerbivory1993}, to name but a few. These HOIs are present in many different systems under different forms, and should therefore not be ignored.

A specific aspect that has been overlooked is the speed at which HOIs affect a system, and how this speed might influence the system dynamics. Consider a system with three species $A$, $B$ and $C$, where the pairwise interaction between $A$ and $B$ is modified by the presence of $C$. We call this HOI a {\em fast} HOI if $C$ modifies the original interaction as soon as it is present, and the modification disappears as soon as species $C$ is removed from the system. On the other hand, we call it a {\em slow} HOI if the modification takes effect gradually, and then in the absence of species $C$ dissipates slowly rather than instantaneously.  
Fast HOIs can be thought of as the manifestation of changing behaviour in the face of competition~\cite{wissingerIntraguildPredationCompetition1993, Letten2019} or direct changes to interaction outcomes~\cite{kumarNaturalHistorydrivenPlantmediated2014, vanveenPlantmodifiedTrophicInteractions2015}.
Slow HOIs can represent, for instance, interaction modifications induced by evolutionary changes~\cite{patelPartitioningEffectsEcoEvolutionary2018, govaertEcoevolutionaryPartitioningMetrics2016} as the third species alters the local environment, or by learning behaviour of certain species \cite{powerWhatCanEcosystems2015}. 
Examples of slow~\cite{jonesPositiveNegativeEffects1997, stoksResurrectingComplexityInterplay2016,sultanResurrectionStudyReveals2013,hendryPerspectivePaceModern1999,manhartGrowthTradeoffsProduce2018} and fast~\cite{ barbosaExperimentalEvidenceHidden2023,liDirectNeighbourhoodEffects2021} HOIs can be found in several ecosystems. However, the effect of the speed of HOIs on the stability of the equilibrium of a system has not yet been studied.

\section*{Mathematical model}
In this paper, HOIs are modelled as an extension of the Generalised Lotka--Volterra model (GLVM), similarly to previous literature~\cite{kleinhesselinkDetectingInterpretingHigherorder2022, singhHigherOrderInteractions2021,aladwaniAdditionHigherorderInteractions2019,terryInteractionModificationsLead2019}. The GLVM is given by the following equations:
\begin{equation}
  \dot{n}_i = n_i \left(1 - n_i + \sum_{j \neq i}^{N} \alpha_{ij} n_j\right)\,,
  \label{GLVM}
\end{equation}
for $i = 1, \ldots , N$, where $n_{i}$ $[-]$~is the abundance of species $i$, $\dot {n}_i$ $[t^{-1}]$~is its first derivative with respect to time, $N$ $[-]$ is the total number of species in the system and $\alpha_{ij}$ $[-]$~is the pairwise interaction coefficient denoting the effect of species~$j$  on species~$i$. The diagonal elements $\alpha_{ii}$ of the interaction matrix are set to zero for convenience, as they represent intraspecific competition, which is already accounted for in the term $-n_i$.

Previous literature~\cite{gibbsCoexistenceDiverseCommunities2022, singhHigherOrderInteractions2021,aladwaniAdditionHigherorderInteractions2019,Letten2019} generally extends the GLVM to include HOIs by adding another term to Eq.~\ref{GLVM}, modelling the multiplicative effect of the abundance of two species. One such extension is
\begin{equation}
  \dot{n}_i = n_i \left(1 - n_i + \sum_{j \neq i}^{N} \alpha_{ij} n_j 
              + \sum_{j\neq i}^{N} \sum_{k\neq i,j}^{N} \beta_{ijk} n_j n_k \right)\,,
  \label{eq: simple HOI}
\end{equation}
where $\beta_{ijk}$ $[-]$ is the HOI coefficient denoting the effect of the combined abundance of species $j$ and $k$ on species $i$ \cite{gibbsCoexistenceDiverseCommunities2022}. In this way, the HOI occurs instantaneously, as soon as species $i$, $j$ and $k$ are present. 

In contrast, here we examine a different kind of HOI, where the HOI does not necessarily occur instantaneously, but its effects can instead emerge gradually. 
This is modelled by introducing a multiplicative modifier $m_{ij}$ $[-]$ to the pairwise interaction $\alpha_{ij}$. The change in the value of $m_{ij}$ is determined by the modifying species $k$. The equations that govern such a system are:
\begin{equation}
  \dot{n}_i  = n_i\left(1-n_i + \sum_{j \neq i}^{N} \alpha_{ij} m_{ij} n_j\right) \label{eq: Natasha ni} \,,
\end{equation}
with, for  $i = 1, \ldots , N$:
\begin{itemize}
\item[(i)] if $\alpha_{ij}$ is modified:
\begin{equation} \label{eq: Natasha mij} 
        \dot{m}_{ij} = \omega \left(1 - m_{ij}  + \sum_{k \neq i,j} \beta_{ijk}  n_k\right)\,;
\end{equation}
\item[(ii)] if $\alpha_{ij}$ is not modified:
\begin{equation}\label{eq: Natasha mij2}
 m_{ij} = 1 \text{ and }  \dot{m}_{ij} = 0\,.
\end{equation}  
\end{itemize}

Here, $\omega$ $[t^{-1} ]$~represents the modification speed, \emph{i.e.}, the speed with which the modifier changes in the presence or absence of the modifying species $k$. $\beta_{ijk}$ is the modification strength. We call $i$ the affected species, $j$ the affecting species and $k$ the modifying species. From an ecological point of view, $\omega$ can be thought of as how quickly the original interactions get changed in the present of the modifying species, while $\beta$ can be thought of as how much the original interactions will get changed as time goes to infinity.
A few observations can be made about the modifier. The first is that the equilibrium value for Eq.~\ref{eq: Natasha mij} is:
\begin{equation}
  \overline{m_{ij}} = 1 + \beta n_k\,,
\end{equation}
which is independent of $\omega$. This means that $\omega$ does not affect the equilibria values of the system, and will only impact the stability of the equilibria and the transient dynamics. If the modifying species $k$ is removed from the system, then the modifier evolves back to $m_{ij} = 1$ and the original pairwise interaction is restored. 

For high values of $\omega$ ($\omega > 100$), the modifier changes much more quickly than the abundances and thus the former can be assumed to be at equilibrium compared to the latter. Therefore, $m_{ij} = 1 + \beta n_k$ for the modified interaction. 
In this case, the equations coincide with those of previously described HOI systems~\cite{gibbsCoexistenceDiverseCommunities2022,singhHigherOrderInteractions2021, Letten2019,aladwaniAdditionHigherorderInteractions2019} if we replace $\beta_{ijk}$ by $\alpha_{ij} \beta_{ijk}$ in Eq.~\ref{eq: simple HOI}:
\begin{equation}
  \dot {n}_i =n_i \left(1 - n_i + \sum_{i \ne j}^{N} \alpha_{ij} n_j + \sum_{j\neq i}^N \sum_{k \ne i,j}^{N} \alpha_{ij} \beta_{ijk} n_i n_j \right)\,.
  \label{eq: simple HOI alpha beta}
\end{equation}

For low values of $\omega$ ($\omega < 0.01$), it can be assumed that the modifier is stationary compared to the species abundances since it the rate of change is much lower. Thus the system simplifies to a GLVM, with the pairwise interaction coefficients modified by the stationary value of the modifier~$\overline{ m_{ij}}$:
\begin{equation}
  \dot{n}_i = n_i \cdot \left(1 - n_i + \sum_{j \neq i}^{N} \overline{ m_{ij}} \alpha_{ij} n_j\right).
  \label{eq: natasha low omega}
\end{equation}

When $\beta = 0$ or $\omega = 0$, there is no modification. In the former case, the modifier will converge to 1, while in the latter case, $m_{ij}$ will be stationary. Since $m_{ij \vert t = 0} = 1$, in both cases $m_{ij} = 1$ at all $t$.
When either of these situations occurs, the system can be described by a classical GLVM as in Eq.~\ref{eq: natasha low omega}, with $\overline{m_{ij}} = 1$. 

\section*{Experimental setting}
\subsection*{Types of HOI interaction networks}
In order to properly study the effect of the modification on a network of interacting species, we limit the scope to networks with only 3 species $A$, $B$ and $C$, with symmetric pairwise interactions. Thus, 
if species $A$ has a negative effect on species $B$, then $B$ has an equally large positive effect on $A$. This means that $\alpha_{ij}  = -\alpha_{ji} $. Furthermore, we assume that all interaction coefficients have the same absolute value $\alpha$: $\vert \alpha_{ij} \vert = \alpha$, for any $i \ne j$.

Restricting to pairwise interactions, two types of networks exist: transitive and intransitive networks. 
In the former, there is a clear hierarchy of competitive superiority: $A > B$, $B>C$ and $A>C$, with $A$ the superior species, $B$ the intermediate and $C$ the inferior one. In the latter case, there is no clear hierarchy: $A > B$, $B>C$ and $C>A$. This means that, if only two out of three species are present, one of them is superior to the other. However, if all three species are present, all species are equal. Substantial research has already been done on these \textit{rock-paper-scissors} systems with pairwise 
interactions~\cite{gallienEffectsIntransitiveCompetition2017,gallienEmergenceWeakintransitiveCompetition2018}, as they can have significant effects on species coexistence~\cite{soliveresEverythingYouAlways2018}. 

\begin{figure}[ht!]
  \centering
  \includegraphics[width=0.45\textwidth]{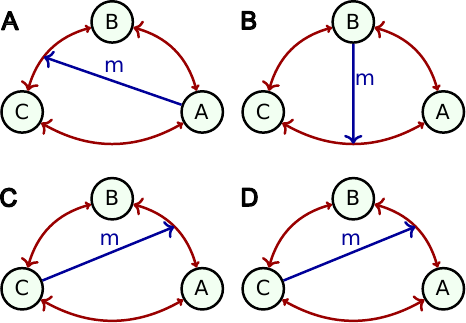}
  \caption{Different 3-species systems with one HOI. The large arrowhead represents the negative side of the interaction, the small arrowhead represents the positive side. (A-C) The three systems where the pairwise interactions form a transitive system. $A > B > C$ and $A > C$ in the pairwise system. (D) A system where the pairwise interactions form an intransitive system. $A > B > C > A$ in the pairwise system.}
  \label{fig: all systems}
\end{figure}

The effect of adding an HOI depends on which pairwise interaction is modified. Since there are three pairwise interactions in the system, three different HOI systems can be distinguished. For the transitive system, these are shown in Figs.~\ref{fig: all systems}A--\ref{fig: all systems}C. For the intransitive system, however, all species are assumed to be equally competitive. This means that modifying any pairwise interaction will have the same effect on the system dynamics. Therefore, only one unique intransitive HOI system exists. This is shown in Fig.~\ref{fig: all systems}D. For the intransitive system, we take $C$ as the modifying species, modifying the interactions between $A$ and $B$.

From now on we will refer to the systems in Fig.~\ref{fig: all systems}A, \ref{fig: all systems}B and \ref{fig: all systems}C as transitive systems $\mathcal{A}$, $\mathcal{B}$ and $\mathcal{C}$, respectively, after the modifying species, and the system in Fig.~\ref{fig: all systems}D as the intransitive system.

\subsection*{Types of HOI}
\begin{figure}[ht!]
  \centering
  \includegraphics[width=0.4\textwidth]{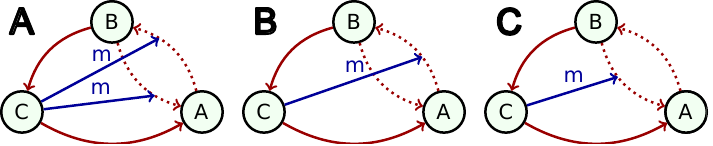}
  \caption{Three different possibilities for one HOI to manifest, shown for the 3-species intransitive system. (A) Symmetric modification~\HOIdouble{ABC}: both sides of the pairwise interaction get modified. (B) Asymmetric modification~\HOI{BAC}: only the negative effect of $A$ on $B$ gets modified. (C) Asymmetric modification~\HOI{ABC}: only the positive effect on $A$ from $B$ gets modified.}
  \label{fig: HOI types}
\end{figure}

Three different types of HOIs are outlined here. Assuming $C$ is the modifying species, the first type is the symmetric HOI, where the modifying species~$C$ modifies the effect of species~$A$ on species~$B$ and vice versa. This is shown in Fig.~\ref{fig: HOI types}A. The two remaining types are asymmetric HOIs, where the modifying species~$C$ only modifies the effect of species $A$ on species $B$, but not the other way around, or $C$ only modifies the effect of species~$B$ on species~$A$. These are shown in Fig.~\ref{fig: HOI types}B and \ref{fig: HOI types}C, respectively. In order to refer to these HOIs more concisely, we use the following notation:
\begin{itemize}
  \item \HOIdouble{ABC}: Both directions of the pairwise interaction between $A$ and $B$ are modified by~$C$ ($\alpha_{AB}$ and $\alpha_{BA}$ are modified).
  \item \HOI{ABC}: Only the effect of species~$B$ on species~$A$ is modified by $C$ ($\alpha_{AB}$ is modified).
  \item \HOI{BAC}: Only the effect of species~$A$ on species~$B$ is modified by $C$ ($\alpha_{BA}$ is modified). 
\end{itemize}

In the remainder of this paper, only the symmetric HOI is discussed unless mentioned otherwise. This is because the symmetric case exhibits interesting dynamics and 
is most appealing from an ecological point of view. 
Note that the asymmetric HOIs are discussed in Supplementary Section~S2.

\section*{Simulations and results}
In this section, we start by examining the effect of the modification speed and strength on the stability of the equilibrium of a system where all pairwise interaction coefficients are identical in absence of modification.  Later on, we drop the restriction that all pairwise interaction coefficients are equal, and look at the effect of the HOI on those systems.

We are specifically interested in the transition between a stable equilibrium and oscillations in the abundance values, for different values the HOI-strength $\beta$ and HOI onset speed $\omega$. The range for $\beta$ is chosen from $0$ to $-80$, as for positive values of $\beta$, no oscillations occur, and values of $\beta \leq  -80$ don't give new results. Negative values for $\beta$ imply the original interaction is weakened ($0 > m > -1$, which implies  $0 > \beta > -9$, see Supplementary Section~S1) or reversed ($-1 > m$, which implies $-9 > \beta$) when $m$ is at equilibrium, and thus the competitively superior species becomes the inferior one and vice versa. While this may seem unrealistic, such competitive reversals have in fact been observed in natural ecosystems~ \cite{chorneskyRepeatedReversalsSpatial1989,steneckMechanismsCompetitiveDominance1991,cantrellCompetitiveReversalsEcological1998,clayEffectsInsectHerbivory1993}. The parameter~$\omega$ is varied between $10^{-3}$ (the modifier is stationary compared to the abundances) and $10^2$ (the modifier is in equilibrium with respect to the abundances).

\subsection*{Identical pairwise interaction coefficients}
None of the three transitive systems shows any oscillations under the addition of an HOI, for any combination of $\alpha$, $\beta$ and $\omega$. This is because when $\alpha > 1$, the competitive pressure on $B$ and $C$ becomes too high, and they go extinct. If $C$ becomes extinct, there is either no pairwise interaction or no interaction modification, depending on whether $C$ is the modifier, and thus no oscillations. When $\alpha < 1$, the pairwise interactions are too weak to induce any oscillations. 

For the intransitive system, we study the effect of parameter combinations of the modification speed~$\omega$ and the modification strength~$\beta$, for a fixed value of $\alpha = 2$. The results for systems with $\alpha = 1$ and $\alpha = 3$ can be found in Supplementary Section~S2.

Figure~\ref{fig: oscillations coexistence intransitive double} shows the effects of $\beta$ and $\omega$ on the stability of the equilibrium and coexistence. In Fig.~\ref{fig: oscillations coexistence intransitive double}A, the blue region shows the parameter combinations for which the system converges to a stable equilibrium point, while the yellow region shows those for which the equilibrium point is unstable. In the latter case, the system converges to a stable limit cycle  and thus shows oscillations. We can see that oscillations only occur for an intermediate range of values of $\omega$. As previously mentioned, for very low values of $\omega$, the system behaves like a GLVM system (Eq.~\ref{GLVM}), and there are no oscillations in the 3-species GLVM system, irrespective of $\beta$. For very high values of $\omega$, the system behaves like an instantaneous HOI model (Eq.~\ref{eq: simple HOI alpha beta}). For this simple HOI model, there are no oscillations for three species. This means that the speed of the HOI is an important factor in the dynamics of the system, and as such cannot be brushed aside.

\begin{figure}[ht!]
  \centering
  \includegraphics[width=0.45\textwidth]{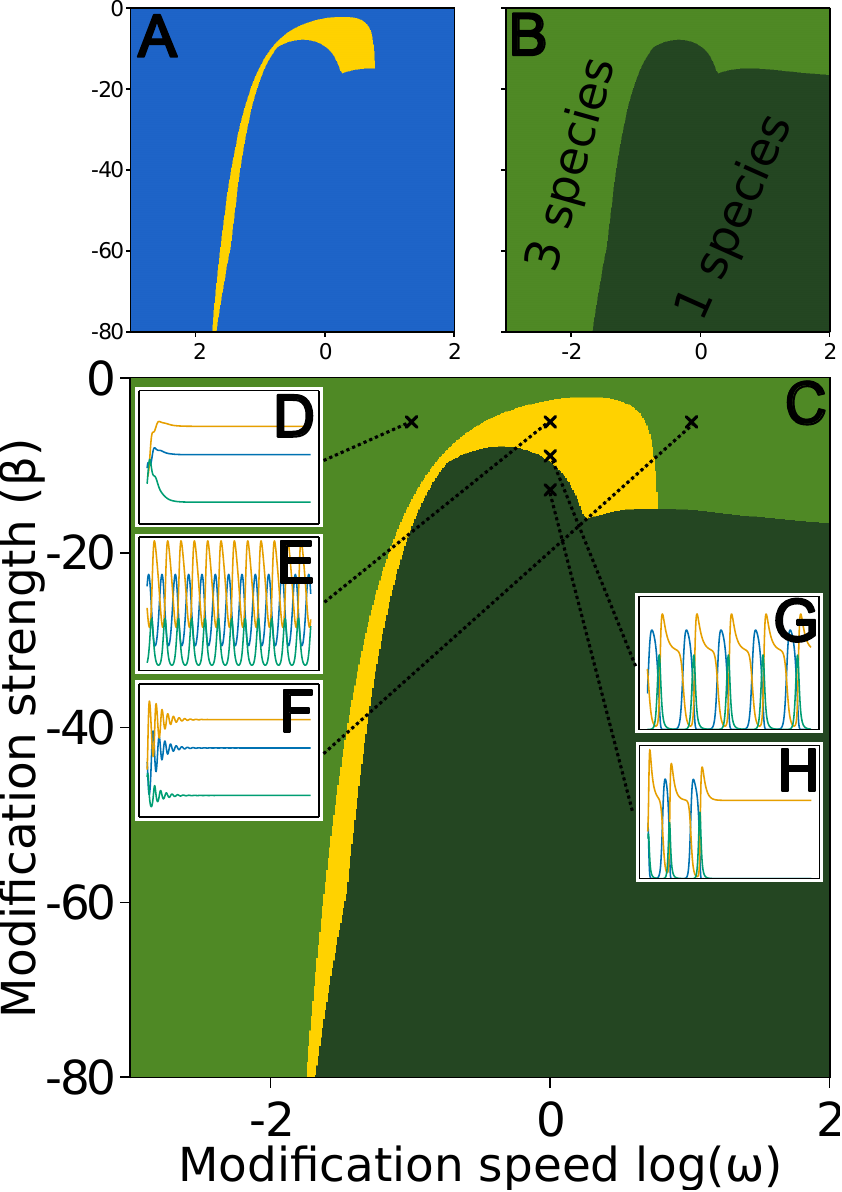}
  \caption{The effect of the modification speed~$\omega$ and modification strength~$\beta$ on the intransitive system stability and coexistence of species. (A) The effect on stability. In the blue region, parameter combinations result in convergence to a stable equilibrium. In the yellow region, parameter combinations result in convergence to a stable limit cycle. (B) The effect on coexistence. The light green region is the region where all species coexist. The dark green region is the region where only one species survives. (C) The superposition of A and B. (D-H) Time series of the species abundances for different parameter combinations. (D) $\omega = 0.1, \beta = -3$ (E) $\omega = 1, \beta = -3$. (F) $\omega = 10, \beta = -3$. (G) $\omega = 1, \beta = -7$. (H) $\omega = 1, \beta = -12$.}
  \label{fig: oscillations coexistence intransitive double}
\end{figure}

Figure~\ref{fig: oscillations coexistence intransitive double}B shows the regions of coexistence. Equilibria with two species do not occur for $\alpha \geq 1$. This is because, if one species becomes extinct, one of the survivors is the superior species and the other the inferior one. When $\alpha \geq 1$, the negative effect of the superior species on the inferior one drives the latter to extinction. In Fig.~\ref{fig: oscillations coexistence intransitive double}C, the region where oscillations occur is superimposed on the coexistence plot. This clearly shows that the oscillations can only be found at the boundary between coexistence and extinction regimes.

As mentioned, $\omega$ has no impact on the equilibrium values of the system. However, for the area to the right of the region of oscillations (\emph{e.g.} Fig.~\ref{fig: oscillations coexistence intransitive double}H) the number of surviving species is always one, while on the left side (\emph{e.g.} Fig.~\ref{fig: oscillations coexistence intransitive double}D) it is three, for the same value of $\beta$. This is because, in a fast system (high $\omega$), one of the species becomes extinct before the system reaches its equilibrium. The fluctuations in the transient abundances are much higher with a higher value of $\omega$. On the left side of the region of oscillations, the abundances start in equilibrium for the fixed value of $m = 1$. As $m$ slowly approaches its equilibrium value, the abundances stay in equilibrium compared to $m$ as $m$ changes much more slowly. This causes the change in abundances during the transient phase to be much slower, which causes the abundances to slowly converge to the equilibrium. The system has more time to reach its equilibrium, and thus the system converges to a stable equilibrium with all species present. The difference between a GLVM (Eq.~\ref{GLVM}) and a simple HOI system (Eq.~\ref{eq: simple HOI}) with the same equilibria becomes clear here. Although they share the same equilibrium, the transient behaviour leads to extinctions in the HOI system, while the GLVM converges to a stable equilibrium with all species present. However, there are cases where the simple HOI system has a stable limit cycle whereas the GLVM has a stable equilibrium. An example of such a system is the intransitive system with HOI \HOI{ABC} with $\alpha > 2$, where oscillations can occur for very high values of $\omega$ (Section~S2). 

Inlaid in Fig.~\ref{fig: oscillations coexistence intransitive double}C are time series of the species abundances for different parameter combinations. For the two examples of stabilising systems without extinctions (D and F), the system converges to the same stable equilibrium, since they only differ in the parameter $\omega$. The two oscillating systems (E and G) have stable limit cycles. When visually comparing these systems, two interesting observations can be made. First, the frequency of the oscillations is much higher for the system with the lower modification strength. Second, the amplitude of the oscillations is higher for the system with the higher modification strength. This can be observed for all parameter combinations: in general, the closer an oscillating system is to the stable region with three coexisting species, the lower the amplitude of the oscillations. The last system shown (H) is an example of a system where the amplitude of the oscillations becomes too large, causing species $A$ to go extinct. Once $A$ goes extinct, the remaining pairwise competition causes species $C$ to go extinct as well since $B > C$. Thus, the closer an oscillating system is to the region with only one species surviving, the higher the amplitude of the oscillations.

\subsection*{Non-identical pairwise interaction coefficients}
Until this point, we have assumed all pairwise interaction coefficients to be equal, and therefore only intransitive systems are prone to oscillations (similar results can be seen for asymmetric HOIs in Section~S2). Transitive systems always converged to a stable equilibrium, for any combination of $\alpha$, $\beta$ and $\omega$. Now, we will look at the effect of non-identical pairwise interaction coefficients on the stability of the equilibrium. In order for oscillations to occur, both the original interaction coefficient and its modification need to be strong enough, without species going extinct. 

\begin{figure*}[tbhp]
  \centering
  \includegraphics[width=.9\textwidth]{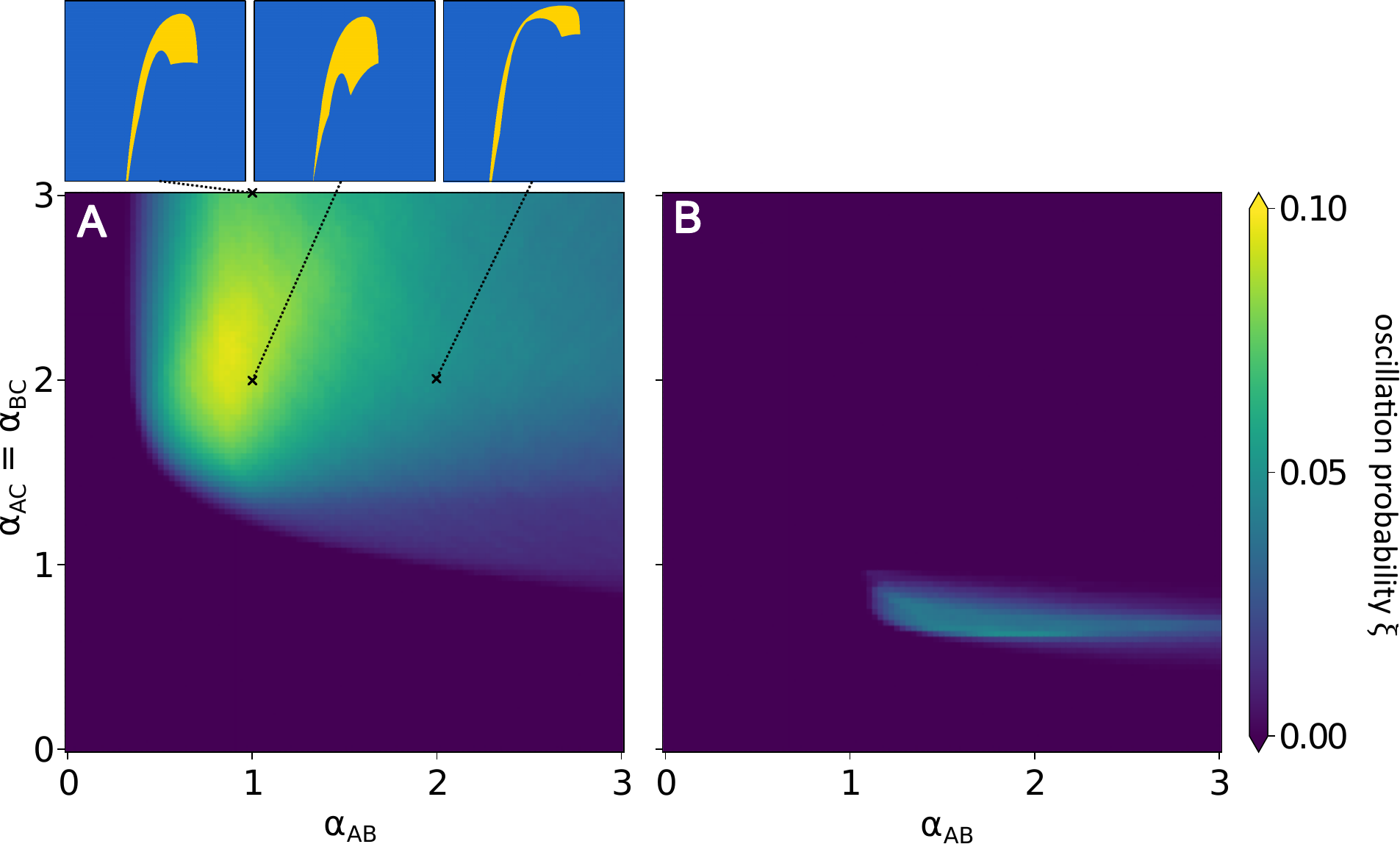}
  \caption{The oscillation probabilities $\xi$ for non-identical pairwise interaction coefficients, where $\alpha_{AB} \neq \alpha_{AC}  = \alpha_{BC} $. (A) The intransitive system. (B) Transitive system $\mathcal{C}$. One pixel corresponds to the oscillation probability $\xi$ for the given combination of $\alpha_{AB}$ and $\alpha_{AC}=\alpha_{BC}$. $\xi$ is defined as the area of oscillations compared to the total area of the parameter space as in Fig.~\ref{fig: oscillations coexistence intransitive double}. The region of oscillations is defined as the area where the system converges to a stable limit cycle. One pixel in this figure captures $27\times 17 = 459$ parameter combinations of $\beta$ and $\omega$, with domains $[10^{-3}, 10^2[$ and $[-80, 0[$ respectively.}
  \label{fig: nonuniform}
\end{figure*}

Therefore, we investigate what happens when $\alpha_{AB}$ is different from $\alpha_{AC} = \alpha_{BC}$ in both the intransitive and transitive system $\mathcal{C}$ (for other cases, see Section~S3). We keep the antisymmetry of the pairwise interactions, so $\alpha_{ij} = -\alpha_{ji}$ for $i,j \in \{A,B,C\}$. However, the pairwise interaction coefficients are not identical. For these systems, we calculated the area of the region in the 
space of the parameters $\beta$ and $\omega$ where oscillations can occur (the yellow area in Fig.~\ref{fig: oscillations coexistence intransitive double}A). Figure~\ref{fig: nonuniform} shows the oscillation probability $\xi$, defined as the area of the region of oscillations, normalised with respect to the total area of the parameter space, for values of $\alpha_{AB}$ and $\alpha_{AC} = \alpha_{BC}$. 

Figure~\ref{fig: nonuniform}A shows the oscillation probabilities for the intransitive system, whereas Fig.~\ref{fig: nonuniform}B shows the oscillation probabilities for transitive system C. 
We can clearly see that the intransitive system is much more prone to oscillation than the transitive system is. The first thing that catches the eye is that for the transitive system, there are no oscillations for $\alpha_{AC} = \alpha_{BC} \geq 1$. This is because in a transitive system, when $\alpha_{AC} \geq 1$, the competitive pressure of $A$ on $C$ is greater than the natural growth rate. Thus, for $\alpha_{AC} > 1$, the equilibrium abundance of $C$ is 0. Since $C$ is the modifier, this means that no modification of the original pairwise interaction can occur, and therefore no oscillations can be introduced by the HOI. Some other transitive systems show no oscillations or are more prone to them, depending on the modifying species and the non-uniform pairwise interaction chosen (Section~S3, Fig.~S2). All intransitive systems have a large area where oscillations occur, for any combination of $\alpha_{AB}$ and $\alpha_{AC}$ (Section~S3, Fig.~S3).

On the other hand, the pairwise interactions need to be strong enough for oscillations to occur. For both systems, we can see that, when all $\vert \alpha_{ij} \vert \leq 1$, there are no oscillations in any of the systems, independent of $\beta$ and $\omega$.

\section*{Discussion}
\subsection*{Causes for (absence of) oscillations}
Oscillations occur for many combinations of values of $\alpha$, $\beta$ and $\omega$ for the intransitive systems. For the intransitive system described above, the underlying mechanism can be understood easily. Consider a point in time when the value of $m$ is high ({\em i.e.}, close to 1, with little modification). Then $A$ experiences a strong positive effect from~$B$, causing $n_{A}$ to increase. In turn, $n_{C}$ increases as well due to the positive effect $A$ has. Lastly, because $n_{C}$ is increasing and $\beta$ is negative, $\dot{m} < 0$ and $m$ starts declining again. The opposite happens once the value of $m$ becomes too low. 

\begin{figure}[h]
  \centering
  \includegraphics[width=.45\textwidth]{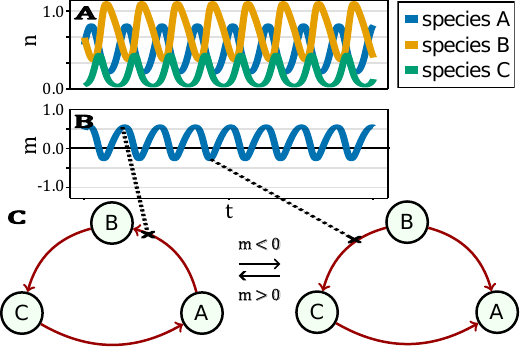}
  \caption{How oscillations occur in the intransitive system ($\beta = -3$, $\omega = 1$). (A) The species abundances $n_i$ over time. (B) The modifier value $m$ over time. (C) When $m < 0$, the system evolves from intransitive to transitive due to modification.}
  \label{fig: system structure change}
\end{figure}

This causes oscillations in the species abundances to occur. From a community point of view, it can be observed that when $m$ is high, the system is intransitive. When $m$ becomes negative, the interaction between $A$ and $B$ reverses, and the system becomes transitive, with $B$ being the superior species. Because of this, $n_{C}$ decreases - causing the absolute value of $m$ to decrease and eventually $m$ becomes positive again. This causes the system to become intransitive again, and the cycle repeats itself. This process is shown in Fig.~\ref{fig: system structure change}.

Oscillations occur in the system when the rates of change of the abundances are sufficiently high. This only happens for combinations of sufficiently high values of $\alpha$ and $\beta$. When these values are too low, the system converges to a stable equilibrium. However, when these values are too high, the amplitude of the oscillations becomes too large, and species go extinct. We can see this in Fig.~\ref{fig: oscillations coexistence intransitive double}, where the region without coexistence lies below the region with oscillations. In the former region, the value of $\beta$ becomes too negative, causing the amplitude of the oscillations to become too large, and species go extinct. 

The speed of the modifier is also important for the oscillations. Although it seems intuitive that higher values of $\omega$ would lead to more oscillations, this is not always the case. We can see this for the intransitive system in Fig.~\ref{fig: oscillations coexistence intransitive double}, where the region of oscillations does not exist for high values of $\omega$. If the value of $\omega$ is too high, then it is possible that $m$ approaches its own equilibrium too quickly, and thus the oscillations disappear. As shown in Fig.~\ref{fig: oscillations coexistence intransitive double}, the oscillations are most pronounced around $\log(\omega) = -0.3$. Here, the amplitude of the oscillations rises very quickly as the absolute value of $\beta$ decreases, resulting in extinctions for values of $\beta$ closer to 0. Thus, in order to for strong oscillations to emerge, the modifier speed should be similar to the natural growth rate of the species. In this case, the delay on the modifier is just long enough to drive oscillations with a higher amplitude.  If the modifier speed is too high, then the delay is too short, and the oscillations disappear. If the modifier speed is too low, then the modifier acts as if it were stationary (additional examples are shown in Section~S2). 

For all of these systems, extinctions occur for less negative values of $\beta$ when $\log(\omega)$ is slightly lower than 1. The effect is more pronounced for  symmetric HOIs compared to asymmetric HOIs. This can be seen when we lower the extinction threshold to $10^{-70}$ instead of $10^{-7}$ (Section~S4). In case of a symmetric HOI, the extinctions still occur for values of $\beta$ close to 0 when $\log(\omega) \approx -0.3$. In case of an asymmetric HOI, the difference in values of $\beta$ for which extinctions occur is not as pronounced, but still present for $\alpha = 2$, and not present for $\alpha = 1$. In the latter case, extinctions due to oscillations no longer occur for $\beta > -70$.

\subsection*{Implications}
While previous theoretical research has shown that HOIs can stabilise ecosystems~\cite{grilliHigherorderInteractionsStabilize2017}, this research suggests that they can be destabilising, especially when the HOI onset speed is taken into account. This follows observations from empirical studies, where HOIs have been shown to cause oscillations in abundances \cite{stoufferCyclicPopulationDynamics2018}.
Oscillations do occur in ecosystems in the real world, independent of yearly phenomena like seasons~\cite{sittlerResponseStoatsMustela1995,ogutuOscillationsLargeMammal2005}. Our models suggest that oscillations of the abundances are closely linked to extinctions, which is important in view of ecosystem conservation, as it can be used to predict the stability of an ecosystem. This is especially essential in the light of the current climate crisis, where ecosystems are under pressure. However, cyclic behaviour can also be beneficial for ecosystem coexistence, as more opportunities for invasion of the other species exist~\cite{stoufferCyclicPopulationDynamics2018}. 
However, it could also be used by humans to destabilise ecosystems. This could be useful in, for example, species management or to counter invasive species. When HOIs are used as a tool for humans to control an ecosystem, it is of vital importance to know the speed of the HOI, which can make the difference between survival and extinction of species.

\subsection*{Limitations and future research}
The results of this paper are limited to a three-species predator-prey network. Future research should be done on how systems with more species change under the influence of HOI onset speed. Especially interesting would be the distinction between even and odd species systems, as these can have vastly different results \cite{gallienEffectsIntransitiveCompetition2017}. Changing the context of the system to a spatially explicit one could also yield interesting results. Here, coexistence may be higher due to migratory behaviour of species. Finally, the results of this paper are limited to predator-prey networks. Future research should be done on how these results change for other types of networks, such as purely competitive or mutualistic networks.

\section*{Conclusion}
Higher order interactions are known to significantly alter the dynamics of species interaction networks. Here, we have demonstrated that the speed at which the HOI manifests itself is a previously overlooked yet important factor in how an ecosystem evolves and on the stability of the equilibrium. We have shown that, after introducing HOI speed as a parameter, for certain speeds of the HOI the system converges to stable oscillations in a three-species system that would otherwise converge to a stable equilibrium. Depending on the 
HOI speed, the system can also converge to an equilibrium with all species present or a single species present. 
By comparing intransitive and transitive systems, we have shown that intransitive systems are more prone to oscillations, as previous research suggests \cite{gilpinLimitCyclesCompetition1975}.
With this research, we have found which parameters and combinations of parameter values are important for a system's equilibrium to become (un)stable. Notably, we have shown that the onset speed of the HOI is one of these parameters that is often vital for (in)stability, and therefore deserves more attention, both in combination with and without intransitive cycles.

\appendix
\section{Methods}
\subsection*{Simulation details}  
Simulations were carried out using Julia~v1.9.3. The code is available on GitHub: \url{https://github.com/Vahieltje/HOIspeednetwork} and Zenodo: \url{https://doi.org/10.5281/zenodo.12706325}. Data/project management in Julia was done with DrWatson~\cite{datserisDrWatsonPerfectSidekick2020}. Plots were made using Makie~\cite{DanischKrumbiegel2021}.
The computational resources (Stevin Supercomputer Infrastructure) and services used in this work were provided by the VSC (Flemish Supercomputer Center), funded by Ghent University, the Research Foundation Flanders and the Flemish Government – Department EWI.

\subsection*{Simulation assumptions}

Ecological systems are initiated with starting abundances $n_i = 1$  since this is the equilibrium abundance for the intransitive system without HOI. The initial value of the modifier was set to $m = 1$, thus assuming no initial modification. The differential equations were solved using Euler's method. The time step is set to $\Delta t = 3 \cdot 10^{-3}$, which is small enough to ensure that the results are not affected by its size. The simulations were run for period that is long enough to ensure the system has reached its equilibrium. The period depends on the modifier speed. For systems with $\omega > 1$, the period was $10\ 000$ time steps. For systems with $\omega<1$, the period is given by $\frac{10\ 000}{\omega}$ time steps, since the modifier needs more time to reach its equilibrium. The system is considered to have reached its equilibrium when:
\begin{equation*}
  \sum_{i=1}^{3} \vert \dot{n}_i \vert < 10^{-4}\,.
\end{equation*}

A species is considered to have gone extinct when its abundance is below the extinction threshold of $n_e = 10^{-7}$. When this happens during the simulation, this species' abundance is set to zero for the remainder of the simulation.

\subsection*{AI tools} To write this paper, T.V.G. used several AI tools. Github-copilot for help when writing the code, Grammarly for help with grammar and spelling, academicgpt to assist writing and spell checking.

\begin{acknowledgments}
  We are grateful to the members of the KERMIT and BionamiX research units and others at the Department of Data Analysis and Mathematical Modelling for useful and interesting discussions, tips and insights. This work was supported by the Research Foundation Flanders under grant number 3G0G0122 and a Ugent-BOF GOA project ``Assessing the biological capacity of ecosystem resilience'',  grant number BOFGOA2017000601.
  \subsection*{Author constributions}
  Conceptualization, All authors; Software, investigation, validation and Writing - original draft, T.V.G.; Funding Acquisition, supervision and writing - Review \& Editing, J.M.B. and B.D.B.; supervision and writing - Review \& Editing, A.J.D. 

  \subsection*{Competing interests}
  There are no competing interests.
\end{acknowledgments}


\bibliography{Citations2}

\begin{thebibliography}{54}%
\makeatletter
\providecommand \@ifxundefined [1]{%
 \@ifx{#1\undefined}
}%
\providecommand \@ifnum [1]{%
 \ifnum #1\expandafter \@firstoftwo
 \else \expandafter \@secondoftwo
 \fi
}%
\providecommand \@ifx [1]{%
 \ifx #1\expandafter \@firstoftwo
 \else \expandafter \@secondoftwo
 \fi
}%
\providecommand \natexlab [1]{#1}%
\providecommand \enquote  [1]{``#1''}%
\providecommand \bibnamefont  [1]{#1}%
\providecommand \bibfnamefont [1]{#1}%
\providecommand \citenamefont [1]{#1}%
\providecommand \href@noop [0]{\@secondoftwo}%
\providecommand \href [0]{\begingroup \@sanitize@url \@href}%
\providecommand \@href[1]{\@@startlink{#1}\@@href}%
\providecommand \@@href[1]{\endgroup#1\@@endlink}%
\providecommand \@sanitize@url [0]{\catcode `\\12\catcode `\$12\catcode
  `\&12\catcode `\#12\catcode `\^12\catcode `\_12\catcode `\%12\relax}%
\providecommand \@@startlink[1]{}%
\providecommand \@@endlink[0]{}%
\providecommand \url  [0]{\begingroup\@sanitize@url \@url }%
\providecommand \@url [1]{\endgroup\@href {#1}{\urlprefix }}%
\providecommand \urlprefix  [0]{URL }%
\providecommand \Eprint [0]{\href }%
\providecommand \doibase [0]{https://doi.org/}%
\providecommand \selectlanguage [0]{\@gobble}%
\providecommand \bibinfo  [0]{\@secondoftwo}%
\providecommand \bibfield  [0]{\@secondoftwo}%
\providecommand \translation [1]{[#1]}%
\providecommand \BibitemOpen [0]{}%
\providecommand \bibitemStop [0]{}%
\providecommand \bibitemNoStop [0]{.\EOS\space}%
\providecommand \EOS [0]{\spacefactor3000\relax}%
\providecommand \BibitemShut  [1]{\csname bibitem#1\endcsname}%
\let\auto@bib@innerbib\@empty
\bibitem [{\citenamefont {{\AA}kesson}\ \emph {et~al.}(2021)\citenamefont
  {{\AA}kesson}, \citenamefont {Curtsdotter}, \citenamefont {Ekl{\"o}f},
  \citenamefont {Ebenman}, \citenamefont {Norberg},\ and\ \citenamefont
  {Barab{\'a}s}}]{akessonImportanceSpeciesInteractions2021}%
  \BibitemOpen
  \bibfield  {author} {\bibinfo {author} {\bibfnamefont {A.}~\bibnamefont
  {{\AA}kesson}}, \bibinfo {author} {\bibfnamefont {A.}~\bibnamefont
  {Curtsdotter}}, \bibinfo {author} {\bibfnamefont {A.}~\bibnamefont
  {Ekl{\"o}f}}, \bibinfo {author} {\bibfnamefont {B.}~\bibnamefont {Ebenman}},
  \bibinfo {author} {\bibfnamefont {J.}~\bibnamefont {Norberg}},\ and\ \bibinfo
  {author} {\bibfnamefont {G.}~\bibnamefont {Barab{\'a}s}},\ }\bibfield
  {title} {\bibinfo {title} {The importance of species interactions in
  eco-evolutionary community dynamics under climate change},\ }\href
  {https://doi.org/10.1038/s41467-021-24977-x} {\bibfield  {journal} {\bibinfo
  {journal} {Nature Communications}\ }\textbf {\bibinfo {volume} {12}},\
  \bibinfo {pages} {4759} (\bibinfo {year} {2021})}\BibitemShut {NoStop}%
\bibitem [{\citenamefont {Chen}\ and\ \citenamefont
  {Cohen}(2001)}]{chenTransientDynamicsFood2001}%
  \BibitemOpen
  \bibfield  {author} {\bibinfo {author} {\bibfnamefont {X.}~\bibnamefont
  {Chen}}\ and\ \bibinfo {author} {\bibfnamefont {J.~E.}\ \bibnamefont
  {Cohen}},\ }\bibfield  {title} {\bibinfo {title} {Transient dynamics and
  food--web complexity in the {{Lotka-Volterra}} cascade model},\ }\href
  {https://doi.org/10.1098/rspb.2001.1596} {\bibfield  {journal} {\bibinfo
  {journal} {Proceedings of the Royal Society of London. Series B: Biological
  Sciences}\ }\textbf {\bibinfo {volume} {268}},\ \bibinfo {pages} {869}
  (\bibinfo {year} {2001})}\BibitemShut {NoStop}%
\bibitem [{\citenamefont {Colon}\ \emph {et~al.}(2015)\citenamefont {Colon},
  \citenamefont {Claessen},\ and\ \citenamefont
  {Ghil}}]{colonBifurcationAnalysisAgentbased2015}%
  \BibitemOpen
  \bibfield  {author} {\bibinfo {author} {\bibfnamefont {C.}~\bibnamefont
  {Colon}}, \bibinfo {author} {\bibfnamefont {D.}~\bibnamefont {Claessen}},\
  and\ \bibinfo {author} {\bibfnamefont {M.}~\bibnamefont {Ghil}},\ }\bibfield
  {title} {\bibinfo {title} {Bifurcation analysis of an agent-based model for
  predator--prey interactions},\ }\href
  {https://doi.org/10.1016/j.ecolmodel.2015.09.004} {\bibfield  {journal}
  {\bibinfo  {journal} {Ecological Modelling}\ }\textbf {\bibinfo {volume}
  {317}},\ \bibinfo {pages} {93} (\bibinfo {year} {2015})}\BibitemShut
  {NoStop}%
\bibitem [{\citenamefont
  {Bunin}(2017)}]{buninEcologicalCommunitiesLotkaVolterra2017}%
  \BibitemOpen
  \bibfield  {author} {\bibinfo {author} {\bibfnamefont {G.}~\bibnamefont
  {Bunin}},\ }\bibfield  {title} {\bibinfo {title} {Ecological communities with
  {{Lotka-Volterra}} dynamics},\ }\href
  {https://doi.org/10.1103/PhysRevE.95.042414} {\bibfield  {journal} {\bibinfo
  {journal} {Physical Review E}\ }\textbf {\bibinfo {volume} {95}},\ \bibinfo
  {pages} {042414} (\bibinfo {year} {2017})}\BibitemShut {NoStop}%
\bibitem [{\citenamefont {Lotka}(1920)}]{lotkaAnalyticalNoteCertain1920}%
  \BibitemOpen
  \bibfield  {author} {\bibinfo {author} {\bibfnamefont {A.~J.}\ \bibnamefont
  {Lotka}},\ }\bibfield  {title} {\bibinfo {title} {Analytical note on certain
  rhythmic relations in organic systems},\ }\href
  {https://doi.org/10.1073/pnas.6.7.410} {\bibfield  {journal} {\bibinfo
  {journal} {Proceedings of the National Academy of Sciences}\ }\textbf
  {\bibinfo {volume} {6}},\ \bibinfo {pages} {410} (\bibinfo {year}
  {1920})}\BibitemShut {NoStop}%
\bibitem [{\citenamefont {Reichenbach}\ \emph {et~al.}(2007)\citenamefont
  {Reichenbach}, \citenamefont {Mobilia},\ and\ \citenamefont
  {Frey}}]{reichenbachMobilityPromotesJeopardizes2007}%
  \BibitemOpen
  \bibfield  {author} {\bibinfo {author} {\bibfnamefont {T.}~\bibnamefont
  {Reichenbach}}, \bibinfo {author} {\bibfnamefont {M.}~\bibnamefont
  {Mobilia}},\ and\ \bibinfo {author} {\bibfnamefont {E.}~\bibnamefont
  {Frey}},\ }\bibfield  {title} {\bibinfo {title} {Mobility promotes and
  jeopardizes biodiversity in rock--paper--scissors games},\ }\href
  {https://doi.org/10.1038/nature06095} {\bibfield  {journal} {\bibinfo
  {journal} {Nature}\ }\textbf {\bibinfo {volume} {448}},\ \bibinfo {pages}
  {1046} (\bibinfo {year} {2007})}\BibitemShut {NoStop}%
\bibitem [{\citenamefont {Suweis}\ \emph {et~al.}(2013)\citenamefont {Suweis},
  \citenamefont {Simini}, \citenamefont {Banavar},\ and\ \citenamefont
  {Maritan}}]{suweisEmergenceStructuralDynamical2013}%
  \BibitemOpen
  \bibfield  {author} {\bibinfo {author} {\bibfnamefont {S.}~\bibnamefont
  {Suweis}}, \bibinfo {author} {\bibfnamefont {F.}~\bibnamefont {Simini}},
  \bibinfo {author} {\bibfnamefont {J.~R.}\ \bibnamefont {Banavar}},\ and\
  \bibinfo {author} {\bibfnamefont {A.}~\bibnamefont {Maritan}},\ }\bibfield
  {title} {\bibinfo {title} {Emergence of structural and dynamical properties
  of ecological mutualistic networks},\ }\href
  {https://doi.org/10.1038/nature12438} {\bibfield  {journal} {\bibinfo
  {journal} {Nature}\ }\textbf {\bibinfo {volume} {500}},\ \bibinfo {pages}
  {449} (\bibinfo {year} {2013})}\BibitemShut {NoStop}%
\bibitem [{\citenamefont {Kessler}\ and\ \citenamefont
  {Shnerb}(2015)}]{kesslerGeneralizedModelIsland2015}%
  \BibitemOpen
  \bibfield  {author} {\bibinfo {author} {\bibfnamefont {D.~A.}\ \bibnamefont
  {Kessler}}\ and\ \bibinfo {author} {\bibfnamefont {N.~M.}\ \bibnamefont
  {Shnerb}},\ }\bibfield  {title} {\bibinfo {title} {Generalized model of
  island biodiversity},\ }\href {https://doi.org/10.1103/PhysRevE.91.042705}
  {\bibfield  {journal} {\bibinfo  {journal} {Physical Review E}\ }\textbf
  {\bibinfo {volume} {91}},\ \bibinfo {pages} {042705} (\bibinfo {year}
  {2015})}\BibitemShut {NoStop}%
\bibitem [{\citenamefont {Arneodo}\ \emph {et~al.}(1982)\citenamefont
  {Arneodo}, \citenamefont {Coullet}, \citenamefont {Peyraud},\ and\
  \citenamefont {Tresser}}]{arneodoStrangeAttractorsVolterra1982}%
  \BibitemOpen
  \bibfield  {author} {\bibinfo {author} {\bibfnamefont {A.}~\bibnamefont
  {Arneodo}}, \bibinfo {author} {\bibfnamefont {P.}~\bibnamefont {Coullet}},
  \bibinfo {author} {\bibfnamefont {J.}~\bibnamefont {Peyraud}},\ and\ \bibinfo
  {author} {\bibfnamefont {C.}~\bibnamefont {Tresser}},\ }\bibfield  {title}
  {\bibinfo {title} {Strange attractors in {{Volterra}} equations for species
  in competition},\ }\href {https://doi.org/10.1007/BF01832841} {\bibfield
  {journal} {\bibinfo  {journal} {Journal of Mathematical Biology}\ }\textbf
  {\bibinfo {volume} {14}},\ \bibinfo {pages} {153} (\bibinfo {year}
  {1982})}\BibitemShut {NoStop}%
\bibitem [{\citenamefont {Vano}\ \emph {et~al.}(2006)\citenamefont {Vano},
  \citenamefont {Wildenberg}, \citenamefont {Anderson}, \citenamefont {Noel},\
  and\ \citenamefont {Sprott}}]{vanoChaosLowdimensionalLotka2006}%
  \BibitemOpen
  \bibfield  {author} {\bibinfo {author} {\bibfnamefont {J.~A.}\ \bibnamefont
  {Vano}}, \bibinfo {author} {\bibfnamefont {J.~C.}\ \bibnamefont
  {Wildenberg}}, \bibinfo {author} {\bibfnamefont {M.~B.}\ \bibnamefont
  {Anderson}}, \bibinfo {author} {\bibfnamefont {J.~K.}\ \bibnamefont {Noel}},\
  and\ \bibinfo {author} {\bibfnamefont {J.~C.}\ \bibnamefont {Sprott}},\
  }\bibfield  {title} {\bibinfo {title} {Chaos in low-dimensional
  {{Lotka}}--{{Volterra}} models of competition},\ }\href
  {https://doi.org/10.1088/0951-7715/19/10/006} {\bibfield  {journal} {\bibinfo
   {journal} {Nonlinearity}\ }\textbf {\bibinfo {volume} {19}},\ \bibinfo
  {pages} {2391} (\bibinfo {year} {2006})}\BibitemShut {NoStop}%
\bibitem [{\citenamefont {Levine}\ \emph {et~al.}(2017)\citenamefont {Levine},
  \citenamefont {Bascompte}, \citenamefont {Adler},\ and\ \citenamefont
  {Allesina}}]{levinePairwiseMechanismsSpecies2017}%
  \BibitemOpen
  \bibfield  {author} {\bibinfo {author} {\bibfnamefont {J.~M.}\ \bibnamefont
  {Levine}}, \bibinfo {author} {\bibfnamefont {J.}~\bibnamefont {Bascompte}},
  \bibinfo {author} {\bibfnamefont {P.~B.}\ \bibnamefont {Adler}},\ and\
  \bibinfo {author} {\bibfnamefont {S.}~\bibnamefont {Allesina}},\ }\bibfield
  {title} {\bibinfo {title} {Beyond pairwise mechanisms of species coexistence
  in complex communities},\ }\href {https://doi.org/10.1038/nature22898}
  {\bibfield  {journal} {\bibinfo  {journal} {Nature}\ }\textbf {\bibinfo
  {volume} {546}},\ \bibinfo {pages} {56} (\bibinfo {year} {2017})}\BibitemShut
  {NoStop}%
\bibitem [{\citenamefont
  {Sanchez}(2019)}]{sanchezDefiningHigherorderInteractions2019}%
  \BibitemOpen
  \bibfield  {author} {\bibinfo {author} {\bibfnamefont {A.}~\bibnamefont
  {Sanchez}},\ }\bibfield  {title} {\bibinfo {title} {Defining higher-order
  interactions in synthetic ecology: Lessons from physics and quantitative
  genetics},\ }\href {https://doi.org/10.1016/j.cels.2019.11.009} {\bibfield
  {journal} {\bibinfo  {journal} {Cell Systems}\ }\textbf {\bibinfo {volume}
  {9}},\ \bibinfo {pages} {519} (\bibinfo {year} {2019})}\BibitemShut {NoStop}%
\bibitem [{\citenamefont {Gibbs}\ \emph {et~al.}(2022)\citenamefont {Gibbs},
  \citenamefont {Levin},\ and\ \citenamefont
  {Levine}}]{gibbsCoexistenceDiverseCommunities2022}%
  \BibitemOpen
  \bibfield  {author} {\bibinfo {author} {\bibfnamefont {T.}~\bibnamefont
  {Gibbs}}, \bibinfo {author} {\bibfnamefont {S.~A.}\ \bibnamefont {Levin}},\
  and\ \bibinfo {author} {\bibfnamefont {J.~M.}\ \bibnamefont {Levine}},\
  }\bibfield  {title} {\bibinfo {title} {Coexistence in diverse communities
  with higher-order interactions},\ }\href
  {https://doi.org/10.1073/pnas.2205063119} {\bibfield  {journal} {\bibinfo
  {journal} {Proceedings of the National Academy of Sciences}\ }\textbf
  {\bibinfo {volume} {119}},\ \bibinfo {pages} {e2205063119} (\bibinfo {year}
  {2022})}\BibitemShut {NoStop}%
\bibitem [{\citenamefont {Terry}\ \emph {et~al.}(2019)\citenamefont {Terry},
  \citenamefont {Morris},\ and\ \citenamefont
  {Bonsall}}]{terryInteractionModificationsLead2019}%
  \BibitemOpen
  \bibfield  {author} {\bibinfo {author} {\bibfnamefont {J.~C.~D.}\
  \bibnamefont {Terry}}, \bibinfo {author} {\bibfnamefont {R.~J.}\ \bibnamefont
  {Morris}},\ and\ \bibinfo {author} {\bibfnamefont {M.~B.}\ \bibnamefont
  {Bonsall}},\ }\bibfield  {title} {\bibinfo {title} {Interaction modifications
  lead to greater robustness than pairwise non-trophic effects in food webs},\
  }\href {https://doi.org/10.1111/1365-2656.13057} {\bibfield  {journal}
  {\bibinfo  {journal} {Journal of Animal Ecology}\ }\textbf {\bibinfo {volume}
  {88}},\ \bibinfo {pages} {1732} (\bibinfo {year} {2019})}\BibitemShut
  {NoStop}%
\bibitem [{\citenamefont {Kleinhesselink}\ \emph {et~al.}(2022)\citenamefont
  {Kleinhesselink}, \citenamefont {Kraft}, \citenamefont {Pacala},\ and\
  \citenamefont {Levine}}]{kleinhesselinkDetectingInterpretingHigherorder2022}%
  \BibitemOpen
  \bibfield  {author} {\bibinfo {author} {\bibfnamefont {A.~R.}\ \bibnamefont
  {Kleinhesselink}}, \bibinfo {author} {\bibfnamefont {N.~J.}\ \bibnamefont
  {Kraft}}, \bibinfo {author} {\bibfnamefont {S.~W.}\ \bibnamefont {Pacala}},\
  and\ \bibinfo {author} {\bibfnamefont {J.~M.}\ \bibnamefont {Levine}},\
  }\bibfield  {title} {\bibinfo {title} {Detecting and interpreting
  higher-order interactions in ecological communities},\ }\href
  {https://doi.org/10.1111/ele.14022} {\bibfield  {journal} {\bibinfo
  {journal} {Ecology Letters}\ }\textbf {\bibinfo {volume} {25}},\ \bibinfo
  {pages} {1604} (\bibinfo {year} {2022})}\BibitemShut {NoStop}%
\bibitem [{\citenamefont {Grilli}\ \emph {et~al.}(2017)\citenamefont {Grilli},
  \citenamefont {Barab{\'a}s}, \citenamefont {{Michalska-Smith}},\ and\
  \citenamefont {Allesina}}]{grilliHigherorderInteractionsStabilize2017}%
  \BibitemOpen
  \bibfield  {author} {\bibinfo {author} {\bibfnamefont {J.}~\bibnamefont
  {Grilli}}, \bibinfo {author} {\bibfnamefont {G.}~\bibnamefont {Barab{\'a}s}},
  \bibinfo {author} {\bibfnamefont {M.~J.}\ \bibnamefont {{Michalska-Smith}}},\
  and\ \bibinfo {author} {\bibfnamefont {S.}~\bibnamefont {Allesina}},\
  }\bibfield  {title} {\bibinfo {title} {Higher-order interactions stabilize
  dynamics in competitive network models},\ }\href
  {https://doi.org/10.1038/nature23273} {\bibfield  {journal} {\bibinfo
  {journal} {Nature}\ }\textbf {\bibinfo {volume} {548}},\ \bibinfo {pages}
  {210} (\bibinfo {year} {2017})}\BibitemShut {NoStop}%
\bibitem [{\citenamefont {Bairey}\ \emph {et~al.}(2016)\citenamefont {Bairey},
  \citenamefont {Kelsic},\ and\ \citenamefont
  {Kishony}}]{baireyHighorderSpeciesInteractions2016a}%
  \BibitemOpen
  \bibfield  {author} {\bibinfo {author} {\bibfnamefont {E.}~\bibnamefont
  {Bairey}}, \bibinfo {author} {\bibfnamefont {E.~D.}\ \bibnamefont {Kelsic}},\
  and\ \bibinfo {author} {\bibfnamefont {R.}~\bibnamefont {Kishony}},\
  }\bibfield  {title} {\bibinfo {title} {High-order species interactions shape
  ecosystem diversity},\ }\href {https://doi.org/10.1038/ncomms12285}
  {\bibfield  {journal} {\bibinfo  {journal} {Nature Communications}\ }\textbf
  {\bibinfo {volume} {7}},\ \bibinfo {pages} {12285} (\bibinfo {year}
  {2016})}\BibitemShut {NoStop}%
\bibitem [{\citenamefont {Xiao}\ \emph {et~al.}(2020)\citenamefont {Xiao},
  \citenamefont {Li}, \citenamefont {Chu}, \citenamefont {Wang}, \citenamefont
  {Meiners},\ and\ \citenamefont
  {Stouffer}}]{xiaoHigherorderInteractionsMitigate2020a}%
  \BibitemOpen
  \bibfield  {author} {\bibinfo {author} {\bibfnamefont {J.}~\bibnamefont
  {Xiao}}, \bibinfo {author} {\bibfnamefont {Y.}~\bibnamefont {Li}}, \bibinfo
  {author} {\bibfnamefont {C.}~\bibnamefont {Chu}}, \bibinfo {author}
  {\bibfnamefont {Y.}~\bibnamefont {Wang}}, \bibinfo {author} {\bibfnamefont
  {S.~J.}\ \bibnamefont {Meiners}},\ and\ \bibinfo {author} {\bibfnamefont
  {D.~B.}\ \bibnamefont {Stouffer}},\ }\bibfield  {title} {\bibinfo {title}
  {Higher-order interactions mitigate direct negative effects on population
  dynamics of herbaceous plants during succession},\ }\href
  {https://doi.org/10.1088/1748-9326/ab8a88} {\bibfield  {journal} {\bibinfo
  {journal} {Environmental Research Letters}\ }\textbf {\bibinfo {volume}
  {15}},\ \bibinfo {pages} {074023} (\bibinfo {year} {2020})}\BibitemShut
  {NoStop}%
\bibitem [{\citenamefont {Raj}\ \emph {et~al.}(2022)\citenamefont {Raj},
  \citenamefont {Upadhyay}, \citenamefont {Karmakar},\ and\ \citenamefont
  {Bhattacharya}}]{rajSimplicialStructuresEcological2022}%
  \BibitemOpen
  \bibfield  {author} {\bibinfo {author} {\bibfnamefont {U.}~\bibnamefont
  {Raj}}, \bibinfo {author} {\bibfnamefont {S.}~\bibnamefont {Upadhyay}},
  \bibinfo {author} {\bibfnamefont {M.}~\bibnamefont {Karmakar}},\ and\
  \bibinfo {author} {\bibfnamefont {S.}~\bibnamefont {Bhattacharya}},\ }\href
  {https://doi.org/10.48550/arXiv.2203.13677} {\bibinfo {title} {Simplicial
  structures in ecological networks}} (\bibinfo {year} {2022}),\ \Eprint
  {https://arxiv.org/abs/2203.13677} {arxiv:2203.13677 [math]} \BibitemShut
  {NoStop}%
\bibitem [{\citenamefont {Terry}\ \emph {et~al.}(2020)\citenamefont {Terry},
  \citenamefont {Bonsall},\ and\ \citenamefont
  {Morris}}]{terryIdentifyingImportantInteraction2020}%
  \BibitemOpen
  \bibfield  {author} {\bibinfo {author} {\bibfnamefont {J.~C.~D.}\
  \bibnamefont {Terry}}, \bibinfo {author} {\bibfnamefont {M.~B.}\ \bibnamefont
  {Bonsall}},\ and\ \bibinfo {author} {\bibfnamefont {R.~J.}\ \bibnamefont
  {Morris}},\ }\bibfield  {title} {\bibinfo {title} {Identifying important
  interaction modifications in ecological systems},\ }\href
  {https://doi.org/10.1111/oik.06353} {\bibfield  {journal} {\bibinfo
  {journal} {Oikos}\ }\textbf {\bibinfo {volume} {129}},\ \bibinfo {pages}
  {147} (\bibinfo {year} {2020})}\BibitemShut {NoStop}%
\bibitem [{\citenamefont {Musciotto}\ \emph {et~al.}(2021)\citenamefont
  {Musciotto}, \citenamefont {Battiston},\ and\ \citenamefont
  {Mantegna}}]{musciottoDetectingInformativeHigherorder2021}%
  \BibitemOpen
  \bibfield  {author} {\bibinfo {author} {\bibfnamefont {F.}~\bibnamefont
  {Musciotto}}, \bibinfo {author} {\bibfnamefont {F.}~\bibnamefont
  {Battiston}},\ and\ \bibinfo {author} {\bibfnamefont {R.~N.}\ \bibnamefont
  {Mantegna}},\ }\bibfield  {title} {\bibinfo {title} {Detecting informative
  higher-order interactions in statistically validated hypergraphs},\ }\href
  {https://doi.org/10.1038/s42005-021-00710-4} {\bibfield  {journal} {\bibinfo
  {journal} {Communications Physics}\ }\textbf {\bibinfo {volume} {4}},\
  \bibinfo {pages} {1} (\bibinfo {year} {2021})}\BibitemShut {NoStop}%
\bibitem [{\citenamefont {Singh}\ and\ \citenamefont
  {Baruah}(2021)}]{singhHigherOrderInteractions2021}%
  \BibitemOpen
  \bibfield  {author} {\bibinfo {author} {\bibfnamefont {P.}~\bibnamefont
  {Singh}}\ and\ \bibinfo {author} {\bibfnamefont {G.}~\bibnamefont {Baruah}},\
  }\bibfield  {title} {\bibinfo {title} {Higher order interactions and species
  coexistence},\ }\href {https://doi.org/10.1007/s12080-020-00481-8} {\bibfield
   {journal} {\bibinfo  {journal} {Theoretical Ecology}\ }\textbf {\bibinfo
  {volume} {14}},\ \bibinfo {pages} {71} (\bibinfo {year} {2021})}\BibitemShut
  {NoStop}%
\bibitem [{\citenamefont {Mayfield}\ and\ \citenamefont
  {Stouffer}(2017)}]{mayfieldHigherorderInteractionsCapture2017}%
  \BibitemOpen
  \bibfield  {author} {\bibinfo {author} {\bibfnamefont {M.~M.}\ \bibnamefont
  {Mayfield}}\ and\ \bibinfo {author} {\bibfnamefont {D.~B.}\ \bibnamefont
  {Stouffer}},\ }\bibfield  {title} {\bibinfo {title} {Higher-order
  interactions capture unexplained complexity in diverse communities},\ }\href
  {https://doi.org/10.1038/s41559-016-0062} {\bibfield  {journal} {\bibinfo
  {journal} {Nature Ecology \& Evolution}\ }\textbf {\bibinfo {volume} {1}},\
  \bibinfo {pages} {1} (\bibinfo {year} {2017})}\BibitemShut {NoStop}%
\bibitem [{\citenamefont {Chatterjee}\ \emph {et~al.}(2022)\citenamefont
  {Chatterjee}, \citenamefont {Nag~Chowdhury}, \citenamefont {Ghosh},\ and\
  \citenamefont {Hens}}]{chatterjeeControllingSpeciesDensities2022}%
  \BibitemOpen
  \bibfield  {author} {\bibinfo {author} {\bibfnamefont {S.}~\bibnamefont
  {Chatterjee}}, \bibinfo {author} {\bibfnamefont {S.}~\bibnamefont
  {Nag~Chowdhury}}, \bibinfo {author} {\bibfnamefont {D.}~\bibnamefont
  {Ghosh}},\ and\ \bibinfo {author} {\bibfnamefont {C.}~\bibnamefont {Hens}},\
  }\bibfield  {title} {\bibinfo {title} {Controlling species densities in
  structurally perturbed intransitive cycles with higher-order interactions},\
  }\href {https://doi.org/10.1063/5.0102599} {\bibfield  {journal} {\bibinfo
  {journal} {Chaos: An Interdisciplinary Journal of Nonlinear Science}\
  }\textbf {\bibinfo {volume} {32}},\ \bibinfo {pages} {103122} (\bibinfo
  {year} {2022})}\BibitemShut {NoStop}%
\bibitem [{\citenamefont {Barbosa}\ \emph {et~al.}(2023)\citenamefont
  {Barbosa}, \citenamefont {Fernandes},\ and\ \citenamefont
  {Morris}}]{barbosaExperimentalEvidenceHidden2023}%
  \BibitemOpen
  \bibfield  {author} {\bibinfo {author} {\bibfnamefont {M.}~\bibnamefont
  {Barbosa}}, \bibinfo {author} {\bibfnamefont {G.~W.}\ \bibnamefont
  {Fernandes}},\ and\ \bibinfo {author} {\bibfnamefont {R.~J.}\ \bibnamefont
  {Morris}},\ }\bibfield  {title} {\bibinfo {title} {Experimental evidence for
  a hidden network of higher-order interactions in a diverse arthropod
  community},\ }\href {https://doi.org/10.1016/j.cub.2022.11.057} {\bibfield
  {journal} {\bibinfo  {journal} {Current Biology}\ }\textbf {\bibinfo {volume}
  {33}},\ \bibinfo {pages} {381} (\bibinfo {year} {2023})}\BibitemShut
  {NoStop}%
\bibitem [{\citenamefont {Li}\ \emph {et~al.}(2021)\citenamefont {Li},
  \citenamefont {Mayfield}, \citenamefont {Wang}, \citenamefont {Xiao},
  \citenamefont {Kral}, \citenamefont {Janik}, \citenamefont {Holik},\ and\
  \citenamefont {Chu}}]{liDirectNeighbourhoodEffects2021}%
  \BibitemOpen
  \bibfield  {author} {\bibinfo {author} {\bibfnamefont {Y.}~\bibnamefont
  {Li}}, \bibinfo {author} {\bibfnamefont {M.~M.}\ \bibnamefont {Mayfield}},
  \bibinfo {author} {\bibfnamefont {B.}~\bibnamefont {Wang}}, \bibinfo {author}
  {\bibfnamefont {J.}~\bibnamefont {Xiao}}, \bibinfo {author} {\bibfnamefont
  {K.}~\bibnamefont {Kral}}, \bibinfo {author} {\bibfnamefont {D.}~\bibnamefont
  {Janik}}, \bibinfo {author} {\bibfnamefont {J.}~\bibnamefont {Holik}},\ and\
  \bibinfo {author} {\bibfnamefont {C.}~\bibnamefont {Chu}},\ }\bibfield
  {title} {\bibinfo {title} {Beyond direct neighbourhood effects: Higher-order
  interactions improve modelling and predicting tree survival and growth},\
  }\href {https://doi.org/10.1093/nsr/nwaa244} {\bibfield  {journal} {\bibinfo
  {journal} {National Science Review}\ }\textbf {\bibinfo {volume} {8}},\
  \bibinfo {pages} {nwaa244} (\bibinfo {year} {2021})}\BibitemShut {NoStop}%
\bibitem [{\citenamefont {{van
  Veen}}(2015)}]{vanveenPlantmodifiedTrophicInteractions2015}%
  \BibitemOpen
  \bibfield  {author} {\bibinfo {author} {\bibfnamefont {F.~F.}\ \bibnamefont
  {{van Veen}}},\ }\bibfield  {title} {\bibinfo {title} {Plant-modified trophic
  interactions},\ }\href {https://doi.org/10.1016/j.cois.2015.02.009}
  {\bibfield  {journal} {\bibinfo  {journal} {Current Opinion in Insect
  Science}\ }\bibinfo {series} {Ecology *
  {{Parasites}}/{{Parasitoids}}/{{Biological}} Control},\ \textbf {\bibinfo
  {volume} {8}},\ \bibinfo {pages} {29} (\bibinfo {year} {2015})}\BibitemShut
  {NoStop}%
\bibitem [{\citenamefont {Morin}\ \emph {et~al.}(1988)\citenamefont {Morin},
  \citenamefont {Lawler},\ and\ \citenamefont
  {Johnson}}]{morinCompetitionAquaticInsects1988}%
  \BibitemOpen
  \bibfield  {author} {\bibinfo {author} {\bibfnamefont {P.~J.}\ \bibnamefont
  {Morin}}, \bibinfo {author} {\bibfnamefont {S.~P.}\ \bibnamefont {Lawler}},\
  and\ \bibinfo {author} {\bibfnamefont {E.~A.}\ \bibnamefont {Johnson}},\
  }\bibfield  {title} {\bibinfo {title} {Competition between aquatic insects
  and vertebrates: Interaction strength and higher order interactions},\ }\href
  {https://doi.org/10.2307/1941637} {\bibfield  {journal} {\bibinfo  {journal}
  {Ecology}\ }\textbf {\bibinfo {volume} {69}},\ \bibinfo {pages} {1401}
  (\bibinfo {year} {1988})}\BibitemShut {NoStop}%
\bibitem [{\citenamefont {Golubski}\ \emph {et~al.}(2016)\citenamefont
  {Golubski}, \citenamefont {Westlund}, \citenamefont {Vandermeer},\ and\
  \citenamefont {Pascual}}]{golubskiEcologicalNetworksEdge2016}%
  \BibitemOpen
  \bibfield  {author} {\bibinfo {author} {\bibfnamefont {A.~J.}\ \bibnamefont
  {Golubski}}, \bibinfo {author} {\bibfnamefont {E.~E.}\ \bibnamefont
  {Westlund}}, \bibinfo {author} {\bibfnamefont {J.}~\bibnamefont
  {Vandermeer}},\ and\ \bibinfo {author} {\bibfnamefont {M.}~\bibnamefont
  {Pascual}},\ }\bibfield  {title} {\bibinfo {title} {Ecological networks over
  the edge: Hypergraph trait-mediated indirect interaction ({{TMII}})
  structure},\ }\href {https://doi.org/10.1016/j.tree.2016.02.006} {\bibfield
  {journal} {\bibinfo  {journal} {Trends in Ecology \& Evolution}\ }\textbf
  {\bibinfo {volume} {31}},\ \bibinfo {pages} {344} (\bibinfo {year}
  {2016})}\BibitemShut {NoStop}%
\bibitem [{\citenamefont {Kumar}\ \emph {et~al.}(2014)\citenamefont {Kumar},
  \citenamefont {Pandit}, \citenamefont {Steppuhn},\ and\ \citenamefont
  {Baldwin}}]{kumarNaturalHistorydrivenPlantmediated2014}%
  \BibitemOpen
  \bibfield  {author} {\bibinfo {author} {\bibfnamefont {P.}~\bibnamefont
  {Kumar}}, \bibinfo {author} {\bibfnamefont {S.~S.}\ \bibnamefont {Pandit}},
  \bibinfo {author} {\bibfnamefont {A.}~\bibnamefont {Steppuhn}},\ and\
  \bibinfo {author} {\bibfnamefont {I.~T.}\ \bibnamefont {Baldwin}},\
  }\bibfield  {title} {\bibinfo {title} {Natural history-driven, plant-mediated
  {{RNAi-based}} study reveals {{CYP6B46}}'s role in a nicotine-mediated
  antipredator herbivore defense},\ }\href
  {https://doi.org/10.1073/pnas.1314848111} {\bibfield  {journal} {\bibinfo
  {journal} {Proceedings of the National Academy of Sciences}\ }\textbf
  {\bibinfo {volume} {111}},\ \bibinfo {pages} {1245} (\bibinfo {year}
  {2014})}\BibitemShut {NoStop}%
\bibitem [{\citenamefont {Wissinger}\ and\ \citenamefont
  {McGrady}(1993)}]{wissingerIntraguildPredationCompetition1993}%
  \BibitemOpen
  \bibfield  {author} {\bibinfo {author} {\bibfnamefont {S.}~\bibnamefont
  {Wissinger}}\ and\ \bibinfo {author} {\bibfnamefont {J.}~\bibnamefont
  {McGrady}},\ }\bibfield  {title} {\bibinfo {title} {Intraguild predation and
  competition between larval dragonflies: Direct and indirect effects on shared
  prey},\ }\href {https://doi.org/10.2307/1939515} {\bibfield  {journal}
  {\bibinfo  {journal} {Ecology}\ }\textbf {\bibinfo {volume} {74}},\ \bibinfo
  {pages} {207} (\bibinfo {year} {1993})}\BibitemShut {NoStop}%
\bibitem [{\citenamefont {Clay}\ \emph {et~al.}(1993)\citenamefont {Clay},
  \citenamefont {Marks},\ and\ \citenamefont
  {Cheplick}}]{clayEffectsInsectHerbivory1993}%
  \BibitemOpen
  \bibfield  {author} {\bibinfo {author} {\bibfnamefont {K.}~\bibnamefont
  {Clay}}, \bibinfo {author} {\bibfnamefont {S.}~\bibnamefont {Marks}},\ and\
  \bibinfo {author} {\bibfnamefont {G.~P.}\ \bibnamefont {Cheplick}},\
  }\bibfield  {title} {\bibinfo {title} {Effects of insect herbivory and fungal
  endophyte infection on competitive interactions among grasses},\ }\href
  {https://doi.org/10.2307/1939935} {\bibfield  {journal} {\bibinfo  {journal}
  {Ecology}\ }\textbf {\bibinfo {volume} {74}},\ \bibinfo {pages} {1767}
  (\bibinfo {year} {1993})}\BibitemShut {NoStop}%
\bibitem [{\citenamefont {Letten}\ and\ \citenamefont
  {Stouffer}(2019)}]{Letten2019}%
  \BibitemOpen
  \bibfield  {author} {\bibinfo {author} {\bibfnamefont {A.~D.}\ \bibnamefont
  {Letten}}\ and\ \bibinfo {author} {\bibfnamefont {D.~B.}\ \bibnamefont
  {Stouffer}},\ }\bibfield  {title} {\bibinfo {title} {The mechanistic basis
  for higher-order interactions and non-additivity in competitive
  communities},\ }\href {https://doi.org/10.1111/ele.13211} {\bibfield
  {journal} {\bibinfo  {journal} {Ecology Letters}\ }\textbf {\bibinfo {volume}
  {22}},\ \bibinfo {pages} {423} (\bibinfo {year} {2019})}\BibitemShut
  {NoStop}%
\bibitem [{\citenamefont {Patel}\ \emph {et~al.}(2018)\citenamefont {Patel},
  \citenamefont {Cortez},\ and\ \citenamefont
  {Schreiber}}]{patelPartitioningEffectsEcoEvolutionary2018}%
  \BibitemOpen
  \bibfield  {author} {\bibinfo {author} {\bibfnamefont {S.}~\bibnamefont
  {Patel}}, \bibinfo {author} {\bibfnamefont {M.~H.}\ \bibnamefont {Cortez}},\
  and\ \bibinfo {author} {\bibfnamefont {S.~J.}\ \bibnamefont {Schreiber}},\
  }\bibfield  {title} {\bibinfo {title} {Partitioning the effects of
  eco-evolutionary feedbacks on community stability},\ }\href
  {https://doi.org/10.1086/695834} {\bibfield  {journal} {\bibinfo  {journal}
  {The American Naturalist}\ }\textbf {\bibinfo {volume} {191}},\ \bibinfo
  {pages} {381} (\bibinfo {year} {2018})}\BibitemShut {NoStop}%
\bibitem [{\citenamefont {Govaert}\ \emph {et~al.}(2016)\citenamefont
  {Govaert}, \citenamefont {Pantel},\ and\ \citenamefont
  {De~Meester}}]{govaertEcoevolutionaryPartitioningMetrics2016}%
  \BibitemOpen
  \bibfield  {author} {\bibinfo {author} {\bibfnamefont {L.}~\bibnamefont
  {Govaert}}, \bibinfo {author} {\bibfnamefont {J.~H.}\ \bibnamefont
  {Pantel}},\ and\ \bibinfo {author} {\bibfnamefont {L.}~\bibnamefont
  {De~Meester}},\ }\bibfield  {title} {\bibinfo {title} {Eco-evolutionary
  partitioning metrics: Assessing the importance of ecological and evolutionary
  contributions to population and community change},\ }\href
  {https://doi.org/10.1111/ele.12632} {\bibfield  {journal} {\bibinfo
  {journal} {Ecology Letters}\ }\textbf {\bibinfo {volume} {19}},\ \bibinfo
  {pages} {839} (\bibinfo {year} {2016})}\BibitemShut {NoStop}%
\bibitem [{\citenamefont {Power}\ \emph {et~al.}(2015)\citenamefont {Power},
  \citenamefont {Watson}, \citenamefont {Szathm{\'a}ry}, \citenamefont {Mills},
  \citenamefont {Powers}, \citenamefont {Doncaster},\ and\ \citenamefont
  {Czapp}}]{powerWhatCanEcosystems2015}%
  \BibitemOpen
  \bibfield  {author} {\bibinfo {author} {\bibfnamefont {D.~A.}\ \bibnamefont
  {Power}}, \bibinfo {author} {\bibfnamefont {R.~A.}\ \bibnamefont {Watson}},
  \bibinfo {author} {\bibfnamefont {E.}~\bibnamefont {Szathm{\'a}ry}}, \bibinfo
  {author} {\bibfnamefont {R.}~\bibnamefont {Mills}}, \bibinfo {author}
  {\bibfnamefont {S.~T.}\ \bibnamefont {Powers}}, \bibinfo {author}
  {\bibfnamefont {C.~P.}\ \bibnamefont {Doncaster}},\ and\ \bibinfo {author}
  {\bibfnamefont {B.}~\bibnamefont {Czapp}},\ }\bibfield  {title} {\bibinfo
  {title} {What can ecosystems learn? {{Expanding}} evolutionary ecology with
  learning theory},\ }\href {https://doi.org/10.1186/s13062-015-0094-1}
  {\bibfield  {journal} {\bibinfo  {journal} {Biology Direct}\ }\textbf
  {\bibinfo {volume} {10}},\ \bibinfo {pages} {69} (\bibinfo {year}
  {2015})}\BibitemShut {NoStop}%
\bibitem [{\citenamefont {Jones}\ \emph {et~al.}(1997)\citenamefont {Jones},
  \citenamefont {Lawton},\ and\ \citenamefont
  {Shachak}}]{jonesPositiveNegativeEffects1997}%
  \BibitemOpen
  \bibfield  {author} {\bibinfo {author} {\bibfnamefont {C.~G.}\ \bibnamefont
  {Jones}}, \bibinfo {author} {\bibfnamefont {J.~H.}\ \bibnamefont {Lawton}},\
  and\ \bibinfo {author} {\bibfnamefont {M.}~\bibnamefont {Shachak}},\
  }\bibfield  {title} {\bibinfo {title} {Positive and negative effects of
  organisms as physical ecosystem engineers},\ }\href
  {https://doi.org/10.1890/0012-9658(1997)078[1946:PANEOO]2.0.CO;2} {\bibfield
  {journal} {\bibinfo  {journal} {Ecology}\ }\textbf {\bibinfo {volume} {78}},\
  \bibinfo {pages} {1946} (\bibinfo {year} {1997})}\BibitemShut {NoStop}%
\bibitem [{\citenamefont {Stoks}\ \emph {et~al.}(2016)\citenamefont {Stoks},
  \citenamefont {Govaert}, \citenamefont {Pauwels}, \citenamefont {Jansen},\
  and\ \citenamefont {De~Meester}}]{stoksResurrectingComplexityInterplay2016}%
  \BibitemOpen
  \bibfield  {author} {\bibinfo {author} {\bibfnamefont {R.}~\bibnamefont
  {Stoks}}, \bibinfo {author} {\bibfnamefont {L.}~\bibnamefont {Govaert}},
  \bibinfo {author} {\bibfnamefont {K.}~\bibnamefont {Pauwels}}, \bibinfo
  {author} {\bibfnamefont {B.}~\bibnamefont {Jansen}},\ and\ \bibinfo {author}
  {\bibfnamefont {L.}~\bibnamefont {De~Meester}},\ }\bibfield  {title}
  {\bibinfo {title} {Resurrecting complexity: The interplay of plasticity and
  rapid evolution in the multiple trait response to strong changes in predation
  pressure in the water flea {{Daphnia}} magna},\ }\href
  {https://doi.org/10.1111/ele.12551} {\bibfield  {journal} {\bibinfo
  {journal} {Ecology Letters}\ }\textbf {\bibinfo {volume} {19}},\ \bibinfo
  {pages} {180} (\bibinfo {year} {2016})}\BibitemShut {NoStop}%
\bibitem [{\citenamefont {Sultan}\ \emph {et~al.}(2013)\citenamefont {Sultan},
  \citenamefont {{Horgan-Kobelski}}, \citenamefont {Nichols}, \citenamefont
  {Riggs},\ and\ \citenamefont {Waples}}]{sultanResurrectionStudyReveals2013}%
  \BibitemOpen
  \bibfield  {author} {\bibinfo {author} {\bibfnamefont {S.~E.}\ \bibnamefont
  {Sultan}}, \bibinfo {author} {\bibfnamefont {T.}~\bibnamefont
  {{Horgan-Kobelski}}}, \bibinfo {author} {\bibfnamefont {L.~M.}\ \bibnamefont
  {Nichols}}, \bibinfo {author} {\bibfnamefont {C.~E.}\ \bibnamefont {Riggs}},\
  and\ \bibinfo {author} {\bibfnamefont {R.~K.}\ \bibnamefont {Waples}},\
  }\bibfield  {title} {\bibinfo {title} {A resurrection study reveals rapid
  adaptive evolution within populations of an invasive plant},\ }\href
  {https://doi.org/10.1111/j.1752-4571.2012.00287.x} {\bibfield  {journal}
  {\bibinfo  {journal} {Evolutionary Applications}\ }\textbf {\bibinfo {volume}
  {6}},\ \bibinfo {pages} {266} (\bibinfo {year} {2013})}\BibitemShut {NoStop}%
\bibitem [{\citenamefont {Hendry}\ and\ \citenamefont
  {Kinnison}(1999)}]{hendryPerspectivePaceModern1999}%
  \BibitemOpen
  \bibfield  {author} {\bibinfo {author} {\bibfnamefont {A.~P.}\ \bibnamefont
  {Hendry}}\ and\ \bibinfo {author} {\bibfnamefont {M.~T.}\ \bibnamefont
  {Kinnison}},\ }\bibfield  {title} {\bibinfo {title} {Perspective: The pace of
  modern life: Measuring rates of contemporary microevolution},\ }\href
  {https://doi.org/10.1111/j.1558-5646.1999.tb04550.x} {\bibfield  {journal}
  {\bibinfo  {journal} {Evolution}\ }\textbf {\bibinfo {volume} {53}},\
  \bibinfo {pages} {1637} (\bibinfo {year} {1999})}\BibitemShut {NoStop}%
\bibitem [{\citenamefont {Manhart}\ and\ \citenamefont
  {Shakhnovich}(2018)}]{manhartGrowthTradeoffsProduce2018}%
  \BibitemOpen
  \bibfield  {author} {\bibinfo {author} {\bibfnamefont {M.}~\bibnamefont
  {Manhart}}\ and\ \bibinfo {author} {\bibfnamefont {E.~I.}\ \bibnamefont
  {Shakhnovich}},\ }\bibfield  {title} {\bibinfo {title} {Growth tradeoffs
  produce complex microbial communities on a single limiting resource},\ }\href
  {https://doi.org/10.1038/s41467-018-05703-6} {\bibfield  {journal} {\bibinfo
  {journal} {Nature Communications}\ }\textbf {\bibinfo {volume} {9}},\
  \bibinfo {pages} {3214} (\bibinfo {year} {2018})}\BibitemShut {NoStop}%
\bibitem [{\citenamefont {AlAdwani}\ and\ \citenamefont
  {Saavedra}(2019)}]{aladwaniAdditionHigherorderInteractions2019}%
  \BibitemOpen
  \bibfield  {author} {\bibinfo {author} {\bibfnamefont {M.}~\bibnamefont
  {AlAdwani}}\ and\ \bibinfo {author} {\bibfnamefont {S.}~\bibnamefont
  {Saavedra}},\ }\bibfield  {title} {\bibinfo {title} {Is the addition of
  higher-order interactions in ecological models increasing the understanding
  of ecological dynamics?},\ }\href {https://doi.org/10.1016/j.mbs.2019.108222}
  {\bibfield  {journal} {\bibinfo  {journal} {Mathematical Biosciences}\
  }\textbf {\bibinfo {volume} {315}},\ \bibinfo {pages} {108222} (\bibinfo
  {year} {2019})}\BibitemShut {NoStop}%
\bibitem [{\citenamefont {Gallien}\ \emph {et~al.}(2017)\citenamefont
  {Gallien}, \citenamefont {Zimmermann}, \citenamefont {Levine},\ and\
  \citenamefont {Adler}}]{gallienEffectsIntransitiveCompetition2017}%
  \BibitemOpen
  \bibfield  {author} {\bibinfo {author} {\bibfnamefont {L.}~\bibnamefont
  {Gallien}}, \bibinfo {author} {\bibfnamefont {N.~E.}\ \bibnamefont
  {Zimmermann}}, \bibinfo {author} {\bibfnamefont {J.~M.}\ \bibnamefont
  {Levine}},\ and\ \bibinfo {author} {\bibfnamefont {P.~B.}\ \bibnamefont
  {Adler}},\ }\bibfield  {title} {\bibinfo {title} {The effects of intransitive
  competition on coexistence},\ }\href {https://doi.org/10.1111/ele.12775}
  {\bibfield  {journal} {\bibinfo  {journal} {Ecology Letters}\ }\textbf
  {\bibinfo {volume} {20}},\ \bibinfo {pages} {791} (\bibinfo {year}
  {2017})}\BibitemShut {NoStop}%
\bibitem [{\citenamefont {Gallien}\ \emph {et~al.}(2018)\citenamefont
  {Gallien}, \citenamefont {Landi}, \citenamefont {Hui},\ and\ \citenamefont
  {Richardson}}]{gallienEmergenceWeakintransitiveCompetition2018}%
  \BibitemOpen
  \bibfield  {author} {\bibinfo {author} {\bibfnamefont {L.}~\bibnamefont
  {Gallien}}, \bibinfo {author} {\bibfnamefont {P.}~\bibnamefont {Landi}},
  \bibinfo {author} {\bibfnamefont {C.}~\bibnamefont {Hui}},\ and\ \bibinfo
  {author} {\bibfnamefont {D.~M.}\ \bibnamefont {Richardson}},\ }\bibfield
  {title} {\bibinfo {title} {Emergence of weak-intransitive competition through
  adaptive diversification and eco-evolutionary feedbacks},\ }\href
  {https://doi.org/10.1111/1365-2745.12961} {\bibfield  {journal} {\bibinfo
  {journal} {Journal of Ecology}\ }\textbf {\bibinfo {volume} {106}},\ \bibinfo
  {pages} {877} (\bibinfo {year} {2018})}\BibitemShut {NoStop}%
\bibitem [{\citenamefont {Soliveres}\ and\ \citenamefont
  {Allan}(2018)}]{soliveresEverythingYouAlways2018}%
  \BibitemOpen
  \bibfield  {author} {\bibinfo {author} {\bibfnamefont {S.}~\bibnamefont
  {Soliveres}}\ and\ \bibinfo {author} {\bibfnamefont {E.}~\bibnamefont
  {Allan}},\ }\bibfield  {title} {\bibinfo {title} {Everything you always
  wanted to know about intransitive competition but were afraid to ask},\
  }\href {https://doi.org/10.1111/1365-2745.12972} {\bibfield  {journal}
  {\bibinfo  {journal} {Journal of Ecology}\ }\textbf {\bibinfo {volume}
  {106}},\ \bibinfo {pages} {807} (\bibinfo {year} {2018})}\BibitemShut
  {NoStop}%
\bibitem [{\citenamefont
  {Chornesky}(1989)}]{chorneskyRepeatedReversalsSpatial1989}%
  \BibitemOpen
  \bibfield  {author} {\bibinfo {author} {\bibfnamefont {E.~A.}\ \bibnamefont
  {Chornesky}},\ }\bibfield  {title} {\bibinfo {title} {Repeated reversals
  during spatial competition between corals},\ }\href
  {https://doi.org/10.2307/1941353} {\bibfield  {journal} {\bibinfo  {journal}
  {Ecology}\ }\textbf {\bibinfo {volume} {70}},\ \bibinfo {pages} {843}
  (\bibinfo {year} {1989})}\BibitemShut {NoStop}%
\bibitem [{\citenamefont {Steneck}\ \emph {et~al.}(1991)\citenamefont
  {Steneck}, \citenamefont {Hacker},\ and\ \citenamefont
  {Dethier}}]{steneckMechanismsCompetitiveDominance1991}%
  \BibitemOpen
  \bibfield  {author} {\bibinfo {author} {\bibfnamefont {R.~S.}\ \bibnamefont
  {Steneck}}, \bibinfo {author} {\bibfnamefont {S.~D.}\ \bibnamefont
  {Hacker}},\ and\ \bibinfo {author} {\bibfnamefont {M.~N.}\ \bibnamefont
  {Dethier}},\ }\bibfield  {title} {\bibinfo {title} {Mechanisms of competitive
  dominance between crustose coralline algae: {{An}} herbivore-mediated
  competitive reversal},\ }\href {https://doi.org/10.2307/1940595} {\bibfield
  {journal} {\bibinfo  {journal} {Ecology}\ }\textbf {\bibinfo {volume} {72}},\
  \bibinfo {pages} {938} (\bibinfo {year} {1991})}\BibitemShut {NoStop}%
\bibitem [{\citenamefont {Cantrell}\ \emph {et~al.}(1998)\citenamefont
  {Cantrell}, \citenamefont {Cosner},\ and\ \citenamefont
  {Fagan}}]{cantrellCompetitiveReversalsEcological1998}%
  \BibitemOpen
  \bibfield  {author} {\bibinfo {author} {\bibfnamefont {R.~S.}\ \bibnamefont
  {Cantrell}}, \bibinfo {author} {\bibfnamefont {C.}~\bibnamefont {Cosner}},\
  and\ \bibinfo {author} {\bibfnamefont {W.~F.}\ \bibnamefont {Fagan}},\
  }\bibfield  {title} {\bibinfo {title} {Competitive reversals inside
  ecological reserves: The role of external habitat degradation},\ }\href
  {https://doi.org/10.1007/s002850050139} {\bibfield  {journal} {\bibinfo
  {journal} {Journal of Mathematical Biology}\ }\textbf {\bibinfo {volume}
  {37}},\ \bibinfo {pages} {491} (\bibinfo {year} {1998})}\BibitemShut
  {NoStop}%
\bibitem [{\citenamefont {Stouffer}\ \emph {et~al.}(2018)\citenamefont
  {Stouffer}, \citenamefont {Wainwright}, \citenamefont {Flanagan},\ and\
  \citenamefont {Mayfield}}]{stoufferCyclicPopulationDynamics2018}%
  \BibitemOpen
  \bibfield  {author} {\bibinfo {author} {\bibfnamefont {D.~B.}\ \bibnamefont
  {Stouffer}}, \bibinfo {author} {\bibfnamefont {C.~E.}\ \bibnamefont
  {Wainwright}}, \bibinfo {author} {\bibfnamefont {T.}~\bibnamefont
  {Flanagan}},\ and\ \bibinfo {author} {\bibfnamefont {M.~M.}\ \bibnamefont
  {Mayfield}},\ }\bibfield  {title} {\bibinfo {title} {Cyclic population
  dynamics and density-dependent intransitivity as pathways to coexistence
  between co-occurring annual plants},\ }\href
  {https://doi.org/10.1111/1365-2745.12960} {\bibfield  {journal} {\bibinfo
  {journal} {Journal of Ecology}\ }\textbf {\bibinfo {volume} {106}},\ \bibinfo
  {pages} {838} (\bibinfo {year} {2018})}\BibitemShut {NoStop}%
\bibitem [{\citenamefont {Sittler}(1995)}]{sittlerResponseStoatsMustela1995}%
  \BibitemOpen
  \bibfield  {author} {\bibinfo {author} {\bibfnamefont {B.}~\bibnamefont
  {Sittler}},\ }\bibfield  {title} {\bibinfo {title} {Response of stoats
  ({{Mustela}} erminea) to a fluctuating lemming ({{Dicrostonyx}}
  groenlandicus) population in {{North East Greenland}}: Preliminary results
  from a long-term study},\ }\href@noop {} {\bibfield  {journal} {\bibinfo
  {journal} {Annales Zoologici Fennici}\ }\textbf {\bibinfo {volume} {32}},\
  \bibinfo {pages} {79} (\bibinfo {year} {1995})},\ \Eprint
  {https://arxiv.org/abs/23735566} {23735566} \BibitemShut {NoStop}%
\bibitem [{\citenamefont {Ogutu}\ and\ \citenamefont
  {{Owen-Smith}}(2005)}]{ogutuOscillationsLargeMammal2005}%
  \BibitemOpen
  \bibfield  {author} {\bibinfo {author} {\bibfnamefont {J.~O.}\ \bibnamefont
  {Ogutu}}\ and\ \bibinfo {author} {\bibfnamefont {N.}~\bibnamefont
  {{Owen-Smith}}},\ }\bibfield  {title} {\bibinfo {title} {Oscillations in
  large mammal populations: Are they related to predation or rainfall?},\
  }\href {https://doi.org/10.1111/j.1365-2028.2005.00587.x} {\bibfield
  {journal} {\bibinfo  {journal} {African Journal of Ecology}\ }\textbf
  {\bibinfo {volume} {43}},\ \bibinfo {pages} {332} (\bibinfo {year}
  {2005})}\BibitemShut {NoStop}%
\bibitem [{\citenamefont {Gilpin}(1975)}]{gilpinLimitCyclesCompetition1975}%
  \BibitemOpen
  \bibfield  {author} {\bibinfo {author} {\bibfnamefont {M.~E.}\ \bibnamefont
  {Gilpin}},\ }\bibfield  {title} {\bibinfo {title} {Limit {{Cycles}} in
  {{Competition Communities}}},\ }\href {https://doi.org/10.1086/282973}
  {\bibfield  {journal} {\bibinfo  {journal} {The American Naturalist}\
  }\textbf {\bibinfo {volume} {109}},\ \bibinfo {pages} {51} (\bibinfo {year}
  {1975})}\BibitemShut {NoStop}%
\bibitem [{\citenamefont {Datseris}\ \emph {et~al.}(2020)\citenamefont
  {Datseris}, \citenamefont {Isensee}, \citenamefont {Pech},\ and\
  \citenamefont {G{\'a}l}}]{datserisDrWatsonPerfectSidekick2020}%
  \BibitemOpen
  \bibfield  {author} {\bibinfo {author} {\bibfnamefont {G.}~\bibnamefont
  {Datseris}}, \bibinfo {author} {\bibfnamefont {J.}~\bibnamefont {Isensee}},
  \bibinfo {author} {\bibfnamefont {S.}~\bibnamefont {Pech}},\ and\ \bibinfo
  {author} {\bibfnamefont {T.}~\bibnamefont {G{\'a}l}},\ }\bibfield  {title}
  {\bibinfo {title} {{{DrWatson}}: The perfect sidekick for your scientific
  inquiries},\ }\href {https://doi.org/10.21105/joss.02673} {\bibfield
  {journal} {\bibinfo  {journal} {Journal of Open Source Software}\ }\textbf
  {\bibinfo {volume} {5}},\ \bibinfo {pages} {2673} (\bibinfo {year}
  {2020})}\BibitemShut {NoStop}%
\bibitem [{\citenamefont {Danisch}\ and\ \citenamefont
  {Krumbiegel}(2021)}]{DanischKrumbiegel2021}%
  \BibitemOpen
  \bibfield  {author} {\bibinfo {author} {\bibfnamefont {S.}~\bibnamefont
  {Danisch}}\ and\ \bibinfo {author} {\bibfnamefont {J.}~\bibnamefont
  {Krumbiegel}},\ }\bibfield  {title} {\bibinfo {title} {Makie.jl: {{Flexible}}
  high-performance data visualization for {{Julia}}},\ }\href
  {https://doi.org/10.21105/joss.03349} {\bibfield  {journal} {\bibinfo
  {journal} {Journal of Open Source Software}\ }\textbf {\bibinfo {volume}
  {6}},\ \bibinfo {pages} {3349} (\bibinfo {year} {2021})}\BibitemShut
  {NoStop}%
\end{thebibliography}%

\end{document}



\title{Supplementary information for Modification speed alters stability of ecological higher-order interaction networks}
\author{Thomas Van Giel, \and Aisling J. Daly \and Jan M. Baetens \and Bernard De Baets}
\date{\today}

\maketitle

\newpage

\section{Analytical equilibria}
\label{sec: analytical equilibria}

\subsection*{Equilibrium solutions}
For some parameter combinations of the three-species system given in Eqs.~3-5 in the main text, the equilibrium values can be calculated analytically. These values do not say anything about the stability of the equilibrium, but they do give an idea of the behaviour of the system. The equilibrium values can be calculated by solving the system's steady state equations, obtained by setting the left-hand side of the differential equations to zero. The equilibrium equations for the intransitive systems are:
\begin{align}
    0 &= n_{A} (1 - n_{A} + \alpha m_1 n_{B} - \alpha n_{C}), \label{eq: nA equilibrium}\\
    0 &= n_{B} (1 - \alpha m_2 n_{A} -  n_{B} + \alpha n_{C}),\label{eq: nB equilibrium}\\
    0 &= n_{C} (1 + \alpha n_{A} - \alpha n_{B} - n_{C}),\label{eq: nC equilibrium}\\
    0 &= \omega (1 - m + \beta n_C),\label{eq: m equilibrium}
\end{align}
where $m_1 = m$ for the systems with HOIs \HOI{ABC} and \HOIdouble{ABC}, $m_1 = 1$ for the system with HOI \HOI{BAC},  $m_2 = m$ for the systems with HOIs \HOI{BAC} and \HOIdouble{ABC}, and $m_2 = 1$ for the system with HOI \HOI{ABC}. 
For $\alpha = 1$, the equilibria for the system with HOI \HOI{ABC} are given by:
\begin{equation}
  n_{A} = \frac{2}{2 - \beta},\quad n_{B} = 1,\quad n_{C} = \frac{2}{2 -\beta},\quad m = \frac{2+\beta}{2-\beta},
  \label{eq: equilibria ABC}
\end{equation}
for $\beta \leq 2$. For the system with HOI \HOI{BAC} the equilibrium is given by:
\begin{equation}
  n_{A} = \frac{2}{2 + \beta},\quad n_{B} = \frac{2}{2 +\beta},\quad n_{C} = 1,\quad m = 1 + \beta,
  \label{eq: equilibria BAC}
\end{equation}
for $\beta \geq -2$. The equilibria in Eq.~\ref{eq: equilibria ABC} and Eq.~\ref{eq: equilibria BAC} can only be reached for $\beta < 2$ and $\beta >  -2$ for \HOI{ABC} and \HOI{BAC}, respectively. Outside of these ranges, the abundances will always grow unbounded. In order to explain this, consider the system with HOI \HOI{ABC}. When $\beta \geq 2$, the positive effect $B$ has on $A$ becomes too strong due to the modification. This causes $n_{A}$ to grow, which in turn causes $n_{C}$ to grow. The growth of $n_{C}$ keeps species $B$ from going extinct due to the competitive pressure of A. Since $n_{C}$ grows, the positive effect $B$ has on $A$ increases as well. This causes $n_{A}$ to grow even faster, which subsequently causes $n_{C}$ to grow even faster, and so on. when $\beta \geq 2$, this feedback loop becomes stronger than the natural death rate so growth becomes unbounded.

In the case of \HOI{BAC}, the negative effect $A$ has on $B$ is reversed when $\beta \leq -1$, and becomes positive. When $\beta \leq -2$, the positive effect $A$ has on $B$ becomes stronger than the natural death rate, while the positive effect of $B$ on $A$ is also stronger than the natural death rate. This implies unbounded growth of $A$ and $B$.

In the symmetric case of \HOIdouble{ABC}, when the negative effect of $A$ onto $B$ is reversed, the same happens to the positive effect of $B$ onto $A$. This means that no positive feedback loop can occur, and the system will never grow indefinitely. Thus, a valid solution exists for all values of $\beta$. 
The general equilibrium equations with variable parameters $\alpha$ and $\beta$ are extensive and not included here due to their length. However, by substituting specific values for $\alpha$ and $\beta$, these equations can be solved.

Because of the unbounded growth of the abundances for the other values of $\beta$, the range over which the analysis is performed is limited to $]-80, 0]$ for the systems with \HOI{ABC} and $[0, 80[$ for the systems with \HOI{BAC}. For the systems with \HOIdouble{ABC}, only systems where $\beta$ lies in $[-80, 0[$ are discussed here. Results from $\beta \in [0, 80[$ are omitted because no oscillations or extinctions occur.

\subsection*{Bifurcation between reversing and weakening the original interaction}
If $\beta < 0$, the modifier either weakens the original interaction or reverses it.  
Finding the bifurcation point between these two states for the intransitive system described in the section \textbf{Identical pairwise interactions} from the main text can be done by substituting $\alpha = 2$ and a value for $m$ in the equilibrium Eqs. \ref{eq: nA equilibrium}-\ref{eq: m equilibrium}. For $1>m>0$, the pairwise interaction is weakened, and $m = 0$ is the bifurcation point where the original interaction is nullified. For $0>m$, the original interaction is reversed. Filling in $m=0$ and $\alpha = 2$, following the system of the section \textbf{Identical pairwise interactions} in the main text, gives the following equilibria for the intransitive system with HOI \HOIdouble{ABC}:

\begin{equation}
  n_{A} = \frac{7}{9},\quad n_{B} = \frac{11}{9},\quad n_{C} = \frac{1}{9},\quad \beta = -9.
\end{equation}

Thus for $0 < \beta < 9$, the original interaction is weakened, whereas for $\beta = 9$ the original interaction is nullified, and for $\beta > 9$ the original interaction is reversed.

\newpage
\section{Nonsymmetric modification of intransitive system with identical values of \texorpdfstring{$\alpha$}{alpha}}
\begin{figure}[h!]
  \centering
  \includegraphics[width = .9\textwidth]{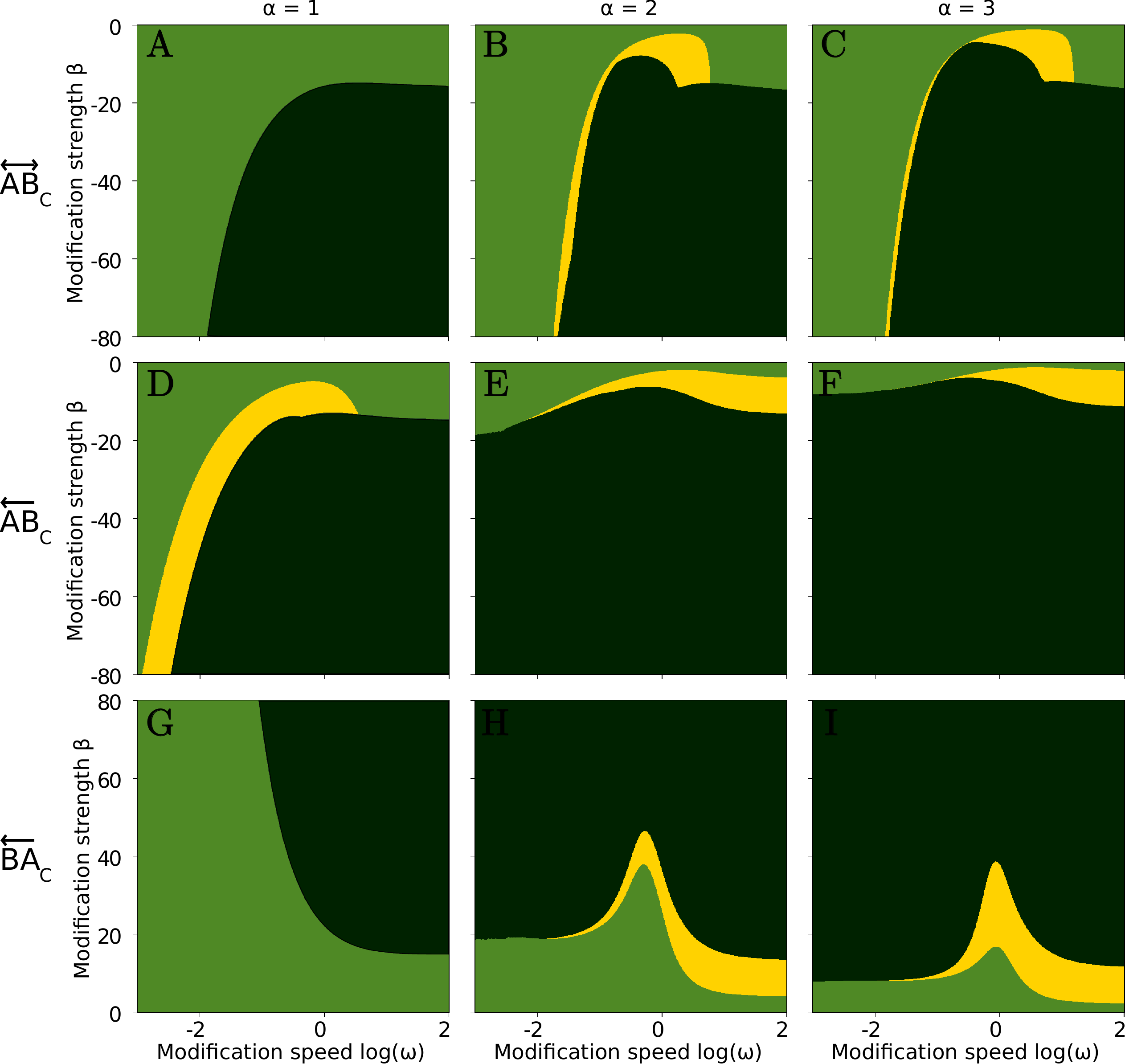}
  \caption{Oscillation regions for symmetric (A-C) and nonsymmetric (D-I) modification, for different values of $\alpha$. The light shade of green corresponds to all three species coexisting, and the dark shade to only one species. The region where oscillations occur are shown in yellow. In this region, all three species coexist as well.}
  \label{fig: oscillations intransitive}
\end{figure}

In the section \textbf{Identical pairwise interactions} of the main text, specific analysis was done for the intransitive system with HOI \HOIdouble{ABC}, with pairwise interaction strength $\alpha = 2$. Here, we highlight the other cases where oscillations occur for systems with identical values of $\alpha$. Since transitive systems with identical values for $\alpha$ do not have oscillations for any combination of $\alpha$, $\beta$ or $\omega$, we limit ourselves to the the intransitive systems. We look at the effect of the modification strength ($\beta$) and speed ($\omega$) for the intransitive system with HOIs \HOIdouble{ABC}, \HOI{ABC} and \HOI{BAC}, and different values of $\alpha \in \{1,2,3\}$ (Fig.~\ref{fig: oscillations intransitive}). 

In this figure, the superposition of the region that allows for oscillations and the coexistence regions are given. The regions where oscillations occur are shown in yellow, and the regions where the system converges to a stable equilibrium are shown in different shades of green, depending on the number of species that coexist. The light shade of green corresponds to the region where all three species coexisting, the dark shade to only one species. In the yellow region with oscillations, all three species coexist as well. The system from the section \textbf{Identical pairwise interactions} from the main text is shown in the Fig.\ref{fig: oscillations intransitive}B. 

The top row in Fig.~\ref{fig: oscillations intransitive} shows the graphs for the intransitive system with symmetric HOI \HOIdouble{ABC}. The system with $\alpha = 1$ does not allow for oscillations. The lowest $\alpha$-value for which oscillations emerge is $\alpha = 1.16$. If $\alpha$ is high enough, oscillations in the systems with a symmetric HOI only occur for an intermediate range of modifier speeds $\omega$. If $\omega$ is either too high or too low, the systems tend towards stable equilibria. For low values of $\omega$ (slow modifier), this equilibrium contains all species. For fast systems, there is a certain threshold for $\beta$ above which the system contains all species, but below which two of the three species go extinct. For systems to the right of the oscillation region ($\log(\omega) > 0.26 \implies \omega > 1.8$), this bifurcation occurs at a relatively constant $\beta$-value, namely $\beta = -16$. For systems below the oscillation region ($\omega < 1.8$), the extinctions start occurring at values of $\beta$ closer to zero. As described in section \textbf{Causes for (absence of) oscillations} in the main text, this \emph{bump} in extinctions for higher values of $\beta$ is because the modifier and species abundances are more synchronised, causing the  oscillation amplitude to increase. The abundance of species $C$ at the lowest point of the oscillation goes below the extinction value, and the species goes extinct. 

The middle row of Fig~\ref{fig: oscillations intransitive} shows the graphs for the intransitive system with asymmetric HOI \HOI{ABC}. Remarkable is that we can see very similar behaviour for the systems with an asymmetric HOI \HOI{ABC} for $\alpha$ equal to 1 as for the systems with a symmetric HOI. In this case (Figure~\ref{fig: oscillations intransitive}D), oscillations also only occur for an intermediate range of values of $\omega$. Contrary to the systems with \HOIdouble{ABC}, the \emph{bump} in values of $\beta$ for which extinctions occur is absent, but at the same values of $\omega$ we can see oscillations occurring for the highest values of $\beta$, pointing to the synchronisation of the modifier and abundances (See also section \ref{sec: changing extinction threshold}). For $\alpha = 2$ and $\alpha = 3$, oscillations occur for a wider range of values of $\omega$, especially for all values of $\omega$ higher than a certain threshold, depending on the value of $\beta$. When $\beta$ gets too negative, extinctions start occurring. For the parameter combinations below the region of oscillations, this is because the oscillation amplitude gets too large, ultimately causing extinctions.  For $\alpha > 2$, the \emph{bump} in values of $\beta$ for which extinctions occur is present, once more indicating the synchronisation between modifier and abundances.

The bottom row of Fig.~\ref{fig: oscillations intransitive} shows the graphs for the intransitive system with asymmetric HOI \HOI{BAC}. As discussed in section~\ref{sec: analytical equilibria}, for these systems the range for $\beta$ is $[0, 80[$ instead of $[-80, 0[$. Although the range of $\beta$ is flipped, we can see that the regions of oscillations and coexistence for systems with \HOI{BAC} is not a mirrored version of the oscillation region for \HOI{ABC} over the $\omega$-axis. This is because the intrinsic growth rate and death rate are not symmetric. Just like in the case with \HOIdouble{ABC}, we do not have oscillations for $\alpha \leq 1$. Once more, the minimum value of $\alpha$ that allows for oscillations is $\alpha = 1.16$. One difference is, for this and higher value of $\alpha$, in the case with \HOIdouble{ABC}, the oscillations occur for intermediate values of $\omega \approx 10^{-0.3}$, while in the case with \HOI{BAC}, we get oscillations for all high values of $\omega > 10^{1.5}$. This shows that the systems with \HOI{BAC} are less stable than the systems with \HOIdouble{ABC}. Since oscillations occur for all high values of $\omega$, the systems with \HOI{BAC} would also show oscillations for the simple HOI-extension, like in Eq.~2 in the main text.
In the systems $\alpha = 2$ or $\alpha= 3$ and HOI \HOI{BAC}, we can see a \emph{bump} in the oscillation region as in the other HOI-systems. However, in the other cases, the \emph{bump} pointed to oscillations occurring for $\beta$ values closer to zero and thus less strong modification is needed to induce oscillations. In this case with \HOI{BAC}, when $\omega \approx 10^{-0.3}$, the values of $\beta$ for which oscillations are possible are much higher than for the surrounding values of $\omega$. Thus, in this case the modification and species abundances being synchronised, prevents oscillations from occurring and prevents species going extinct instead of facilitating them.

\newpage
\section{Non-identical pairwise interactions}

Here we look at the effect of non-identical pairwise interaction strengths $\alpha_{ij}$ on the oscillation probabilities for all cases, both with transitive and intransitive systems. 
For easier reference, let's call a system with $\alpha_{AB} \neq \alpha_{AC}  = \alpha_{BC} $ system $\widehat{AB}$, since $\alpha_{AB}$ is the pairwise interaction with a different value compared to the other two. 
Pairwise interactions are still symmetric, so that $\alpha_{ij} = \alpha_{ji}$. 

Multiple different parameter combinations exist when working with non-identical $\alpha$ values. The different cases are:
\begin{itemize}
  \item Transitive system A,B or C, or intransitive (4 possibilities)
  \item symmetric or either of two asymmetric HOI (3 possibilities)
  \item $\widehat{AB}$, $\widehat{AC}$ or $\widehat{BC}$. (3 possibilities)
\end{itemize}

This comes down to 36 different combinations, most of which do not result in any oscillations for irrespective of parameter combination. In the main text, we only looked at transitive system A and the intransitive system, with a symmetric HOI and $\widehat{AB}$. Table~\ref{tab: oscillating systems} shows the combinations for which oscillations occur.
\begin{table}[ht]
  \centering
  \caption{Do any oscillating systems exist for different modifiers, non-identical pairwise interactions and HOIs?}
  \label{tab: oscillating systems}
  \begin{tabular}{lllll}
                                           &                                & \multicolumn{3}{c}{\shortstack{non-identical\\ interaction}}                                     \\ \cline{3-5} 
  \rule{0pt}{3ex}                           & \multicolumn{1}{c|}{HOI-type}  & \multicolumn{1}{c|}{$\widehat{AB}$} & \multicolumn{1}{c|}{$\widehat{AC}$} & \multicolumn{1}{l|}{$\widehat{BC}$} \\ \cline{2-5} 
  \multicolumn{1}{l|}{Transitive system A} & \multicolumn{1}{l|}{\rule{0pt}{3ex} \HOIdouble{BCA}}    & No                      & No                      & No                      \\ \cline{2-2}
  \multicolumn{1}{c|}{}                    & \multicolumn{1}{l|}{\rule{0pt}{3ex} \HOI{BCA}} & No                      & No                      & No                      \\ \cline{2-2}
  \multicolumn{1}{l|}{}                    & \multicolumn{1}{l|}{\rule{0pt}{3ex} \HOI{CBA}} & No                      & No                      & No                      \\ \cline{2-2}
                                           &                                &                         &                         &                         \\ \cline{2-2}
  \multicolumn{1}{c|}{Transitive system B} & \multicolumn{1}{l|}{\rule{0pt}{3ex} \HOIdouble{ACB}}    & No                      & No                      & No                      \\ \cline{2-2}
  \multicolumn{1}{l|}{}                    & \multicolumn{1}{l|}{\rule{0pt}{3ex} \HOI{ACB}} & No             & Yes                     & No                      \\ \cline{2-2}
  \multicolumn{1}{l|}{}                    & \multicolumn{1}{l|}{\rule{0pt}{3ex} \HOI{CAB}} & No                      & No                      & No                      \\ \cline{2-2}
                                           &                                &                         &                         &                         \\ \cline{2-2}
  \multicolumn{1}{l|}{Transitive system C} & \multicolumn{1}{l|}{\rule{0pt}{3ex} \HOIdouble{ABC}}    & Yes                     & Yes                     & No                      \\ \cline{2-2}
  \multicolumn{1}{l|}{}                    & \multicolumn{1}{l|}{\rule{0pt}{3ex} \HOI{ABC}} & No                      & Yes                     & No                      \\ \cline{2-2}
  \multicolumn{1}{c|}{}                    & \multicolumn{1}{l|}{\rule{0pt}{3ex} \HOI{BAC}} & No                      & No                      & No                      \\ \cline{2-2}
  \multicolumn{1}{c}{}                     &                                &                         &                         &                         \\ \cline{2-2}
  \multicolumn{1}{l|}{Intransitive system} & \multicolumn{1}{l|}{\rule{0pt}{3ex} \HOIdouble{ABC}}    & Yes                     & Yes                     & Yes                     \\ \cline{2-2}
  \multicolumn{1}{l|}{}                    & \multicolumn{1}{l|}{\rule{0pt}{3ex} \HOI{ABC}} & Yes                     & Yes                     & Yes                     \\ \cline{2-2}
  \multicolumn{1}{l|}{}                    & \multicolumn{1}{l|}{\rule{0pt}{3ex} \HOI{BAC}} & Yes                     & Yes                     & Yes                     \\ \cline{2-2}
  \end{tabular}
  \end{table}

\subsection*{Transitive systems}

For the transitive systems, oscillations do not occur for system A. In this system, $A$ modifies the interaction between $B$ and $C$. In the transitive system, $B$ and $C$ are the two least abundant species, which means that their pairwise interaction has the least effect on the system. This is why the modification of this interaction cannot induce any oscillations.
On the other hand, transitive system C is the transitive system with that allows for oscillations most often. Here the interaction between the two most abundant species is modified. This means that the modification has the most effect on the system. 

\begin{figure}[ht!]
  \centering
  \includegraphics[width=0.6\textwidth]{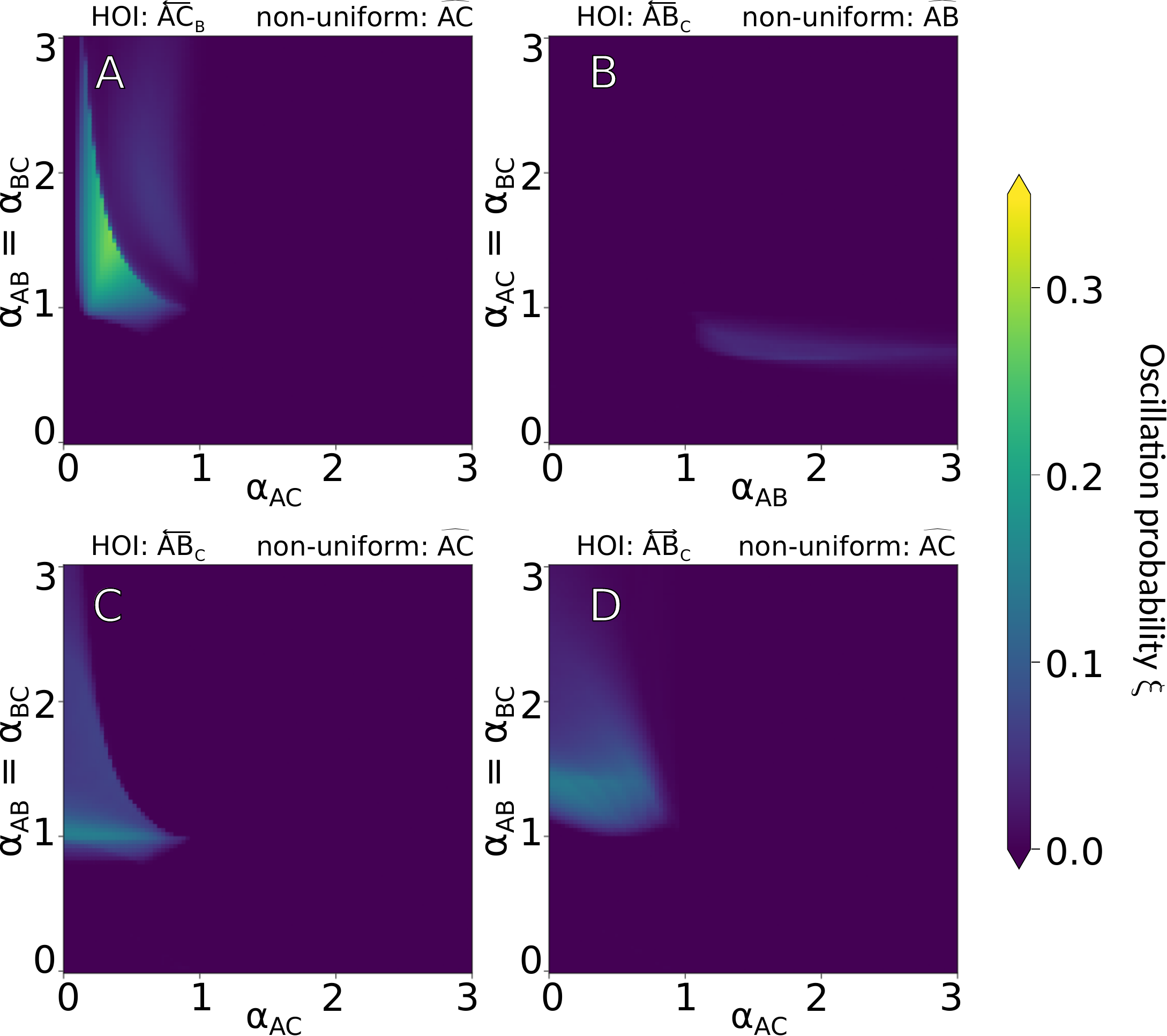}
  \caption{Oscillation probabilities $\xi$ for all transitive systems listed in Table~\ref{tab: oscillating systems} with non-identical $\alpha$ values, where oscillations occur.}
  \label{fig: nonidentical alphas transitive}
\end{figure}

The oscillation probabilities $\xi$ are shown in Fig.~\ref{fig: nonidentical alphas transitive}. We can see that oscillations only occur for systems where an interaction influencing the abundance of $A$ is modified. This is because $A$ is the most abundant species in the transitive systems, and thus has the most influence on how the system evolves. The modifications of all these systems reduce or invert the positive effects of other species onto $A$. This means that the abundance of $A$ will drop, and the abundance of the other species will rise. The abundances of the other species are then high enough to influence the abundance of $A$ again, providing them with the possibilities of inducing oscillations.

Just like described in section \textbf{Non-identical pairwise interactions} of the main text, we can see in Figure~\ref{fig: nonidentical alphas transitive} that for the transitive systems oscillations only occur for $\alpha_{AC} < 1$. If $\alpha_{AC} > 1$, the competitive pressure of $A$ on $C$ becomes to high, leading $C$ to go extinct, and thus being unable to induce oscillations. 

From the four transitive systems allowing oscillations, three of them have non-identical pairwise interaction $\widehat{AC}$, and one has $\widehat{AB}$. For the systems with $\widehat{AC}$, the oscillation probability $\xi$ is much higher. The difference between $\widehat{AC}$ and $\widehat{AB}$ is that, in the systems with $\widehat{AC}$, when $\alpha_{AC}$ is low, $\alpha_{BC}$ can still be high. We can see that, for $\widehat{AC}$, $\xi$ is higher for higher values of $\alpha_{BC}$. Thus the interaction between $B$ and $C$ is a key driver in the oscillations of the system.

\subsection*{Intransitive systems}

\begin{figure}[ht!]
  \centering
  \includegraphics[width=0.8\textwidth]{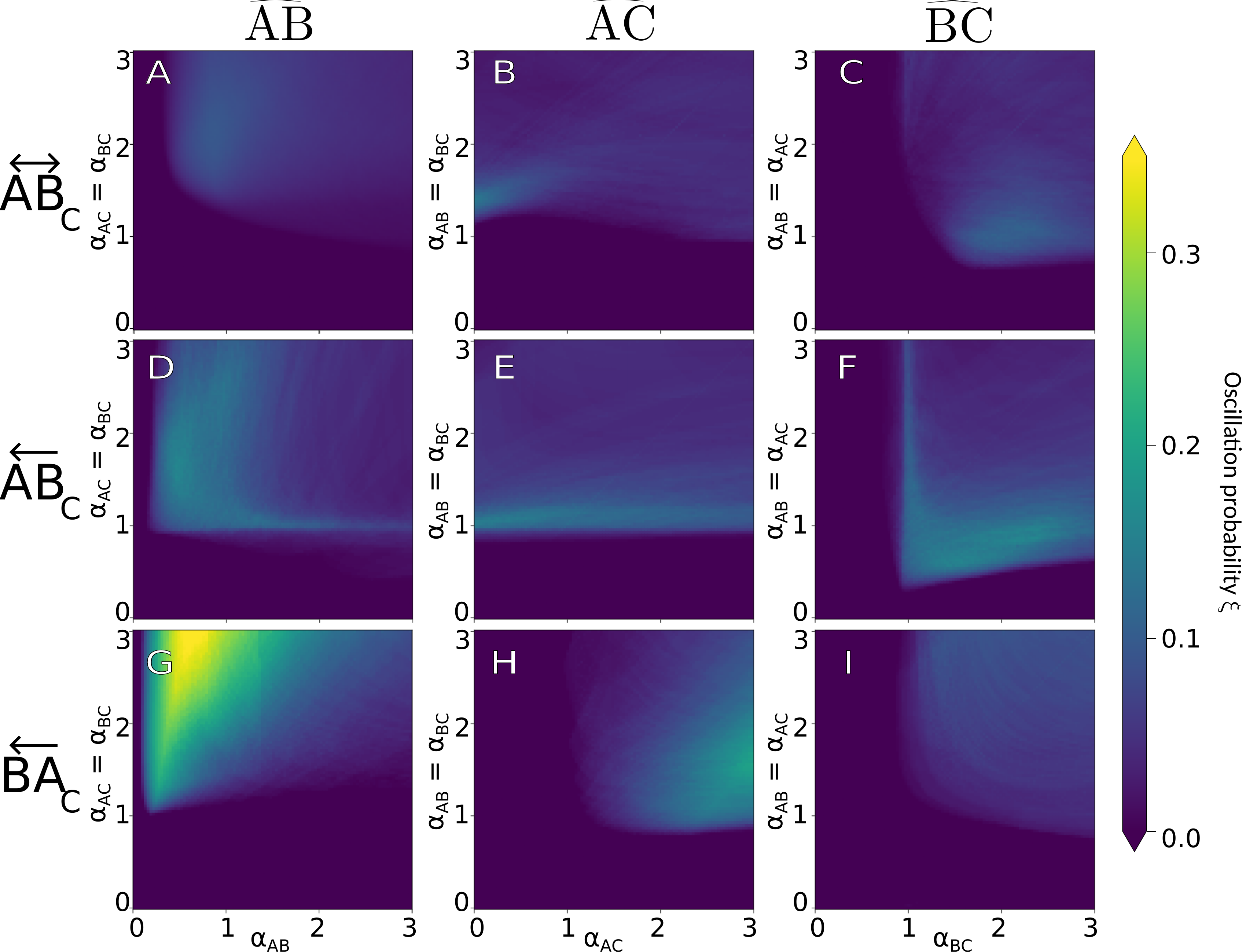}
  \caption{Oscillation probabilities $\xi$ for all intransitive systems listed in Table~\ref{tab: oscillating systems}.}
  \label{fig: nonuniform alphas intransitive}
\end{figure}

As we can see from Table~\ref{tab: oscillating systems} any combination of HOI-type and non-identical interaction can cause oscillations to occur in the intransitive system. 
The oscillation probabilities $\xi$ for the intransitive systems are shown in Figure \ref{fig: nonuniform alphas intransitive}.

For the intransitive systems, oscillations can occur for a wide range of values of $\alpha$ and HOI-types.
We can see that the ranges of parameter combinations for the different pairwise interaction strengths for which oscillations occur are much bigger than for the transitive systems. We can see that, for all intransitive systems, when all different interaction strengths $\alpha_{ij}$ are high enough, oscillations occur.

The first thing that can be observed is that systems with HOI \HOIdouble{ABC} and \HOI{ABC} have similar behaviour on how the oscillation probability $\xi$ changes with the respect to $\alpha$, especially when compared to \HOI{BAC}. 

For these systems, multiple interesting observations can be made. The first is that $\alpha_{AC}$ does not seem to affect the oscillation probability $\xi$ much. Indeed, when looking at systems with $\widehat{AC}$ (Figs.~\ref{fig: nonuniform alphas intransitive}B and E), we see that $\xi$ is quite constant for different values of $\alpha_{AC}$. We can also see this effect when comparing the systems with $\widehat{AB}$ and $\widehat{BC}$ in Figures~\ref{fig: nonuniform alphas intransitive}A and C or in Figures~\ref{fig: nonuniform alphas intransitive}D and F. These are very similar figures upon switching the axes. The most important pairwise interaction seems to be $\alpha_{AB}$ as the change in $\xi$ is more pronounced for changes in $\alpha_{AB}$ compared to $\alpha_{AC}$ and $\alpha_{BC}$. We can see that $\xi$ always first increases with $\alpha_{AB}$, and then decreases. This is because, when $\alpha_{AB}$ is high, extinctions start to happen as the oscillation amplitude gets too large.

Lastly, we can see that systems with \HOI{BAC} have much higher values for $\xi$, compared to the other systems, especially for $\widehat{AB}$ and $\widehat{AC}$. In these systems, when the values for $\alpha$ increases, the oscillation amplitude does not increase as much, so extinctions do not occur as rapidly as in systems with other HOIs. This means that the oscillations are much more stable for high values of $\alpha$ compared to intermediate values.


\newpage
\section{Changing the extinction threshold}
\label{sec: changing extinction threshold}
\begin{figure}[ht!]
  \centering
  \includegraphics[width=0.8\textwidth]{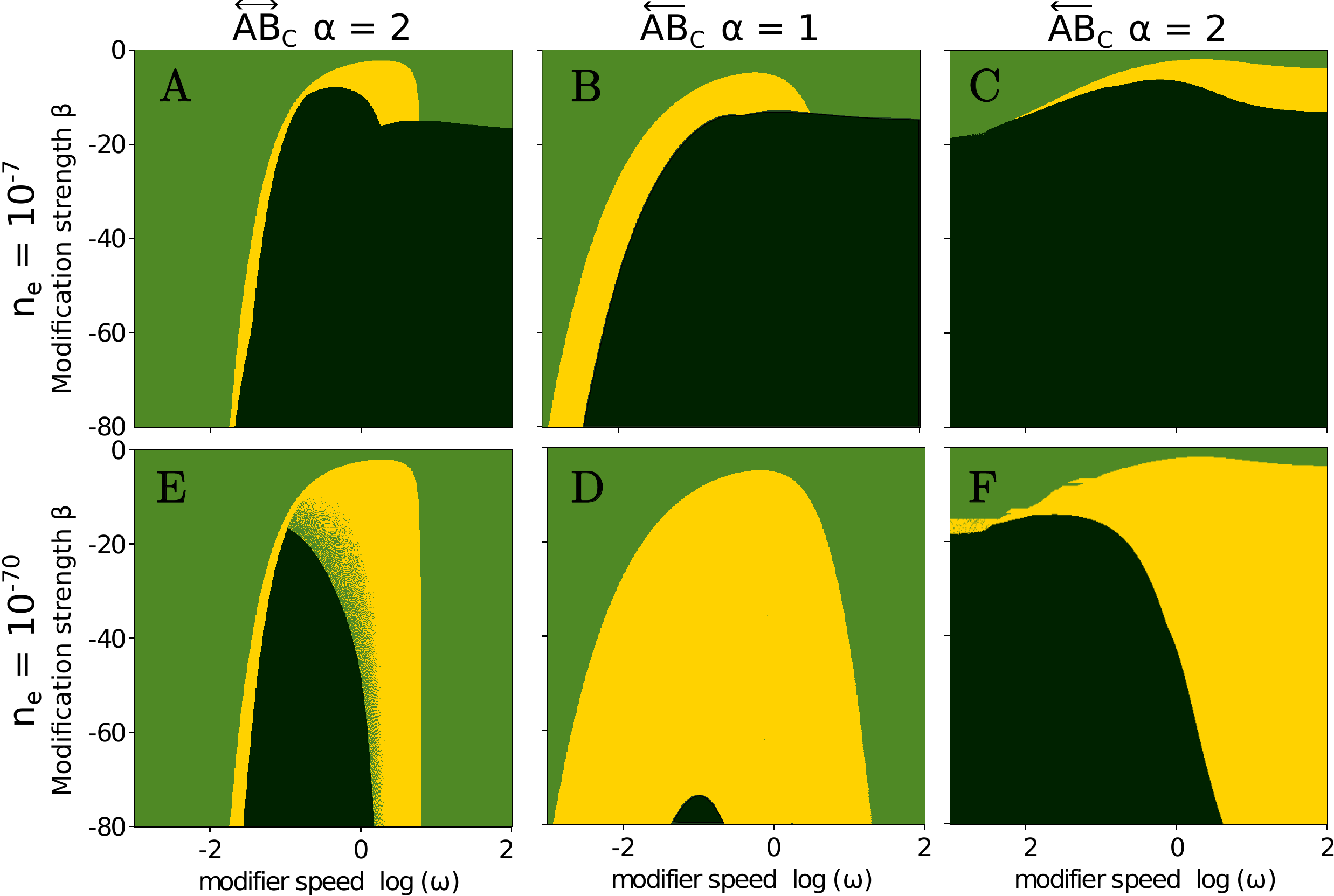}
  \caption{The effect of changing the extinction threshold $n_e$ of intransitive systems. The light green region is the region where no extinctions occur. The dark green region is the region where 2 species become extinct. The yellow region is the region where all species coexist, but the equilibrium is a limit cycle instead of a fixed point.-}
  \label{fig: changing extinction threshold}
\end{figure}

In the main text and the previous sections of the supplementary material, we have assumed that the extinction threshold is equal to $n_e = 10^{-7}$, which means that, if a species abundance $n_{i} \leq n_e$, the species $i$ is considered extinct and its abundance is set to $0$. In this section, we show the effects of the extinction threshold on the behaviour of the systems. Figure~\ref{fig: changing extinction threshold} shows the oscillation maps of some of the intransitive systems from section S2, with different thresholds for extinction.

The top row in figure~\ref{fig: changing extinction threshold} are the oscillation maps of systems where the extinction threshold equals $n_e = 10^{-7} $. The bottom row shows the oscillation maps of same systems where the extinction threshold is set to $10^{-70}$. 
For the systems with HOI \HOI{ABC} and $\alpha = 1$ (Figure~\ref{fig: changing extinction threshold}E), we can see that the region of oscillations extends almost all the way to the bottom of the figure. Thus the value of $\beta$ needs to be very negative to induce oscillations with high enough amplitudes to cause extinctions. 

On the other hand, for systems with \HOIdouble{ABC} (Figure \ref{fig: changing extinction threshold}D) we can still see extinctions occurring below the oscillation regions for much higher values of $\beta$, meaning that the amplitudes rise much faster with $\beta$, compared to the other systems. Especially around $\log(\omega) = -1$, the amplitude of the oscillations seems to rise very quickly with $\beta$. 

For systems with \HOI{ABC} and $\alpha = 2$ (Figure \ref{fig: changing extinction threshold}F), we can see that the oscillation region extends all the way to the bottom of the figure for high values of $\omega$, but for intermediate values ($-1.7<\log(\omega)<0)$ we can again see that extinctions happen for relatively high values of $\beta$. This is again because the oscillation amplitude gets too high, causing species to go extinct. Thus in this kind of system, the oscillations are much more stable (i.e. the oscillation amplitude does not explode with $\beta$) for high values of $\omega$ compared compared to moderate values.